\documentclass[sigconf]{acmart}
\copyrightyear{2019} 
\acmYear{2019} 
\setcopyright{acmcopyright}
\acmConference[ISSAC '19] {2019 International Symposium on Symbolic and Algebraic Computation}{July 15--18, 2019}{Beijing, China}

	\RequirePackage{calc} 
	\RequirePackage{enumerate} 
	\RequirePackage{framed}	
	\newcommand{\Xchee}[1]{}
	\newcommand{\Xmichael}[1]{}
	\newcommand{\Xjuan}[1]{}
\newcommand{\newOld}[2]{#1}
\newcommand{\ifNeedSpace}[1]{}		

	\newcommand{\btable}[2][ccccccccccccccccccccccccc]{
	    \begin{center}
		\begin{tabular}{#1}
		#2
		\end{tabular}
	    \end{center}}
	
	\newcommand{\refPro}[1]{Proposition~\protect{\ref{pro:#1}}}
	
	\newcommand{\0}{\mbox{\bf 0}}  	
	\newcommand{\1}{\mbox{\bf 1}}  	
	\newcommand{\true}{\mbox{\tt true}}	
	\newcommand{\false}{\mbox{\tt false}}	
		
	\newcommand{\half}{\textstyle{\frac{1}{2}}} 
		\newcommand{\aast}{^{\ast}}	
	\newcommand{\INPUt}{\mbox{\textsf{Input}}}
	\newcommand{\OUTPUt}{\mbox{\textsf{Output}}}
	\newcommand{\argmin}{\mathop{\textrm{argmin}}}
	\definecolor{shadecolor}{RGB}{220,255,220}	
	\newcommand{\bshade}{\begin{shaded}}	
	\newcommand{\eshade}{\end{shaded}}	
	\newcommand{\zero}{\mbox{\texttt{Zero}}}     

		\newcommand{\mtf}{\ensuremath{\texttt{MT}_\bff}}
	\newcommand{\mkf}{\ensuremath{\texttt{MK}_\bff}}	
		\newcommand{\mk}{{\ensuremath{\texttt{MK}}}}	
		\newcommand{\intmk}{{\ensuremath{\intbox\texttt{MK}}}}
		\newcommand{\wtmk}{\ensuremath{\wtbox{\texttt{MK}}}}
	\newcommand{\jc}{{\texttt{JC}}}			
    \newcommand{\jcs}{{\texttt{JC}_\textrm{s}}}		
		\newcommand{\intjc}{{\ensuremath{\intbox\texttt{JC}}}}

		\newcommand{\wtjc}{\ensuremath{\wtbox{\texttt{JC}}}}

	\newcommand{\cz}{{\ensuremath{\texttt{C}_0}}}	
		\newcommand{\intcz}{{\ensuremath{\intbox\texttt{C}_0}}}
		\newcommand{\wtcz}{\ensuremath{\wtbox{\texttt{C}_0}}}
    \newcommand{\K}[1]{K(#1)}
	
	\newcommand{\sep}{{\mathrm sep}}
		\newcommand{\boxm}{\intbox_{_{\rm M}}}	
        \newcommand{\wtboxm}{\wt{\intbox}_{_{\rm M}}}

    \newcommand{\jkd}[2]{\|J\inv(\Delta(#1))\cdot K(\Delta(#2))\|}
    \newcommand{\discalpha}{\Delta_\bfalpha}

	\newcommand{\Comment}[1]{
	    \quad\textcolor{cyan}{\mbox{$\triangleleft\;\;$}{\em #1}} }
	\newcommand{\bfa}{{\boldsymbol{a}}}

	\newcommand{\bfalpha}{{\boldsymbol{\alpha}}}
	\newcommand{\bfb}{{\boldsymbol{b}}}

	\newcommand{\bff}{{\boldsymbol{f}}}

	\newcommand{\bfg}{{\boldsymbol{g}}}

	\newcommand{\bfm}{{\boldsymbol{m}}}
	
	\newcommand{\bfp}{{\boldsymbol{p}}}

	\newcommand{\bfx}{{\boldsymbol{x}}}

	\newcommand{\bfy}{{\boldsymbol{y}}}

	\newcommand{\bfz}{{\boldsymbol{z}}}

		\newcommand{\calS}{{\mathcal S}}

		\newcommand{\calZ}{{\mathcal Z}}

	\newcommand{\RR}{{\mathbb R}}
	
	\newcommand{\ZZ}{{\mathbb Z}}
	\newcommand{\as}{{\mathrel{\,:=\,}}}
	\newcommand{\dd}{ ,\ldots , }
	\newcommand{\wt}[1]{\widetilde{#1}}
	\newcommand{\wtO}{\widetilde{O}}	
	\newcommand{\wh}[1]{\widehat{#1}}
	\newcommand{\ib}{\subseteq }
	\newcommand{\ignore}[1]{}
	
	\newcommand{\ii}{\boldsymbol{i}}	
	\newcommand{\dt}[1]{\textbf{#1}}
	\newcommand{\ass}{\leftarrow}	
	\newcommand{\inv}{^{-1}}	
	\newcommand{\es}{\emptyset}
	\newcommand{\beq}{\begin{equation}}
	\newcommand{\eeq}{\end{equation}}
	\newcommand{\beql}[1]{\begin{equation}\label{eq:#1}}
	\newcommand{\eeql}{\end{equation}}
	\newcommand{\beqarray}{\begin{eqnarray}}
	\newcommand{\eeqarray}{\end{eqnarray}}	
	\newcommand{\beqarrays}{\begin{eqnarray*}}
	\newcommand{\eeqarrays}{\end{eqnarray*}} 
	\newcommand{\set}[1]{\left\{ #1 \right\}}
	\newcommand{\eps}{\varepsilon}
	\newcommand{\vareps}{\varepsilon}
	\newcommand{\efrac}[1]{\frac{1}{#1}}	

	\newcommand{\beqarrayl}[1]{\begin{eqnarray}\label{eq:#1}}
	\newcommand{\eeqarrayl}{\end{eqnarray}}
	\newcommand{\blem}{\begin{mylemma}}
	\newcommand{\elem}{\end{mylemma}}
	\newcommand{\bleml}[1]{\begin{mylemma} \label{lem:#1}}
	\newcommand{\eleml}{\end{mylemma}}
	\newcommand{\blemT}[2]{\begin{mylemma}[#1] \label{lem:#2}}
	\newcommand{\elemT}{\end{mylemma}}
	\newcommand{\bthmT}[2]{\begin{mytheorem}[#1] \label{thm:#2}}
	\newcommand{\ethmT}{\end{mytheorem}}
	\newcommand{\refeq}[1]{(\protect{\ref{eq:#1}})}

	\newcommand{\refLem}[1]{Lemma~\protect{\ref{lem:#1}}}
	\newcommand{\refThm}[1]{Theorem~\protect{\ref{thm:#1}}}

    \newtheorem{mytheorem}{\textsc{Theorem}}[section]
	\newtheorem{mylemma}[mytheorem]{\textsc{Lemma}}
	\newtheorem{mycoro}[mytheorem]{\textsc{Corollary}}
	\newtheorem{mypropo}[mytheorem]{\textsc{Proposition}}
	\newcommand{\bthml}[1]{\begin{mytheorem} \label{thm:#1}}
	\newcommand{\ethml}{\end{mytheorem}}
	\newcommand{\bpro}{\begin{mypropo}}
	\newcommand{\epro}{\end{mypropo}}
	\newcommand{\bprol}[1]{\begin{mypropo} \label{pro:#1}}
	\newcommand{\eprol}{\end{mypropo}}
	\newcommand{\bproT}[2]{\begin{mypropo}[#1] \label{pro:#2}}
	\newcommand{\eproT}{\end{mypropo}}
	\newcommand{\bcor}{\begin{mycoro}}
	\newcommand{\ecor}{\end{mycoro}}
	\newcommand{\bcorl}[1]{\begin{mycoro} \label{cor:#1}}
	\newcommand{\ecorl}{\end{mycoro}}
	\newcommand{\bcorT}[2]{\begin{mycoro}[#1] \label{cor:#2}}
	\newcommand{\ecorT}{\end{mycoro}}
 	\newcommand{\bgenDIY}[3][0.15in]{
 		\vspace*{#1}
 		\noindent \textbf{#2}
 		\textit{ #3 }\ \\
 		}
	\newenvironment{progb}[2][4]{ 
	\begin{center}
	\fbox{\begin{minipage}{0.75\textwidth}
	\vspace*{-#1\abovedisplayskip}
	\begin{prog}#2\end{prog}
	\end{minipage}}
	\end{center}
	}{}
	\newcommand{\myprogb}[1]{\begin{progb}{ #1 }\end{progb}}

	\newenvironment{prog}{\begin{tabbing}
xxxx\=xx\=xx\=xx\=xx\=xx\=xxxx\=xxxx\=xxxx\=xxxx\=xxxx\=xxxx\=xxxx\=
	\kill\\}{
	\end{tabbing}}
	\newcommand{\lline}[1][0]{
	\ \\ \Indent[#1]\texttt {}
	}
	
	\newcounter{indentcounter1}
	\newcounter{indentcounter2}
	\newcounter{indentcounter3}
	\newcommand{\Indent}[1][5]{
	\setcounter{indentcounter1}{#1}		
	\setcounter{indentcounter2}{1}		
	\setcounter{indentcounter3}{\value{indentcounter1}*
			\value{indentcounter2}}
	\hspace*{\value{indentcounter3} mm}}
	\newcommand{\bpf}[1][Proof.]{\begin{pf}[#1]~}
	\newcommand{\epf}{\end{pf}}
	\newenvironment{pf}[1][Proof.]{{\em #1}}{ 
	\hspace*{1mm}\hfill \textbf{ Q.E.D.} \vspace{2mm} \noindent}
			
	\newsavebox{\doubleveesym}
	\sbox{\doubleveesym}{
	        \begin{picture}(8,6)(1,1)
	                \put(1,1){$\vee$}
	                \put(3,1){$\vee$}
	        \end{picture}
	        }
	
	\newsavebox{\doublewedgesym}
	\sbox{\doublewedgesym}{
	        \begin{picture}(8,6)(1,1)
	                \put(1,1){$\wedge$}
	                \put(3,1){$\wedge$}
	        \end{picture}
	        }
	
	\newsavebox{\doublebigveesym}
	\savebox{\doublebigveesym}{
	        \begin{picture}(12,9)(1,2)
	                \put(1,4){$\bigvee$}
	                \put(4,4){$\bigvee$}
	        \end{picture}
	        }
	
	\newsavebox{\doublebigwedgesym}
	\savebox{\doublebigwedgesym}{
	        \begin{picture}(12,9)(1,2)
	                \put(1,4){$\bigwedge$}
	                \put(4,4){$\bigwedge$}
	        \end{picture}
	        }
	\newcommand{\doublebigwedge}{\mathop{\usebox{\doublebigwedgesym}}}
	\newcommand{\grouping}[2][lllllllllllllllllllllllll]{\left.
	\begin{array}{#1}
	#2
	\end{array}\right\}}
	\newsavebox{\intboxsym}
	\newcommand{\intbox}{
	    {\,\,\setlength{\unitlength}{.33mm}\framebox(4,7){}\,}}
	\newcommand{\wtbox}{{\wt{\intbox}}}
	\newcommand{\intboxm}{\intbox_{_{\tt M}}}	
	\newcommand{\bitem}{\begin{itemize}}
	\newcommand{\eitem}{\end{itemize}}
	\newcommand{\benum}{\begin{enumerate}}
	\newcommand{\eenum}{\end{enumerate}}
	\newcommand{\bdescr}{\begin{description}}
	\newcommand{\edescr}{\end{description}}
 
 \newcommand{\miranda}{{\tt Miranda}}	

\title[Isolating Zeros of Real Systems of Equations]{
	Effective Subdivision Algorithm for Isolating Zeros\\
	of Real Systems of Equations, with Complexity Analysis}

\author{Juan Xu}
\authornote{
    Partially supported by the National Natural Science Foundation of 
    China (NSFC 11771034).}
	\affiliation{%
	\institution{Beihang University}
	\streetaddress{37, Xueyuan Road}
	\city{Beijing}
	\state{China}
	\postcode{100191}}
	\email{xujuan@buaa.edu.cn}
\author{Chee Yap}
	\authornote{Partially supported by NSF Grants
	Nos.~CCF-1423228 and CCF-1564132.
	Further support under Chinese Academy of Science (Beijing)
	President's International Fellowship Initiative (2018),
	and Beihang International Visiting Professor Program No. Z2018060.}
	\affiliation{%
	\institution{Courant Institute, NYU}
	\streetaddress{251 Mercer Street}
	\city{New York}
	\state{NY}
	\postcode{10012}}
	\email{yap@cs.nyu.edu}
\date{\today}

\begin{abstract}
    We describe a new algorithm \miranda\ for isolating the simple
    zeros of a function $\bff:\RR^n\to\RR^n$
    within a box $B_0\ib \RR^n$.  The function $\bff$ and its partial
    derivatives must have interval forms, but need not be polynomial.
    Our subdivision-based algorithm is ``effective'' in the
    sense that our algorithmic description also specifies the
    numerical precision that is sufficient to certify an implementation
    with any standard BigFloat number type. 
    The main predicate is the Moore-Kioustelides (MK) test,
    based on Miranda's Theorem (1940).
    Although the MK test is well-known, this paper appears to be the first
    synthesis of this test into a complete root isolation algorithm.
	
    We provide a complexity analysis of
    our algorithm based on intrinsic geometric parameters
    of the system. 
    Our algorithm and complexity analysis are 
    developed using 3 levels of
    description (Abstract, Interval, Effective).
    This methodology provides a systematic pathway for
    achieving effective subdivision algorithms in general. 
\end{abstract}

\keywords{Root Isolation; 
	System of Real Equations;
	Certified Computation; 
	Subdivision Algorithms;
	Miranda Theorem;
	Effective Certified Algorithm;
	Complexity Analysis;}

\begin{document}

\maketitle

\section{Introduction}
	Solving multivariate zero-dimensional systems of equations is a fundamental
	task with many applications.  We focus on the problem of
	isolating simple real zeros of a real function
		$$\bff=(f_1\dd f_n):\RR^n\to\RR^n$$
	within a given bounded box $B_0\ib\RR^n$.
	We do not require $\bff$ to be polynomial, only each
	$f_i$ and its partial derivatives have interval forms.
	We require that $\bff$ has only isolated simple zeros
    in $2B_0$\newOld{, which is the box sharing the same center 
    of $B_0$ and of with twice that of $B_0$}{}.
	We call $B_0$ the region-of-interest (ROI) of the input instance.
	This formulation of root isolation is called\footnote{
		Sometimes, an algorithm is called ``local'' if it works
		in small enough neighborhoods (like Newton iteration),
		and ``global'' if no such restriction is needed.
		Clearly, this is a different local/global distinction.
	}
	a \dt{local problem}
	in \cite{imbach-pan-yap:ccluster:18},
	in contrast to the \dt{global problem} of
	isolating all roots of $\bff$.
	The local problem is very
	important in higher dimensions because the global problem has
	complexity that is exponential in $n$.  In geometric applications we
	typically can identify ROI's and can solve the corresponding local
	problem much faster than the global problem.  Moreover, if $\bff$ is
	not polynomial, the global problem might not be solvable:
	E.g., $\bff = \sin x$, $n=1$.  But it is solvable as
	a local problem as in \cite{yap-sagraloff-sharma:cluster:13}.
	
	In their survey of root finding in polynomial systems,
	Sherbrooke and Patrikalakis \cite{sherbrooke-patrikalakis:93}
	noted 3 main approaches:
	(1) algebraic techniques, (2) homotopy, (3) subdivision. 
	They objected to the first two approaches on ``philosophical grounds'',
	meaning that it is not easy in these methods to restrict its
	computation to some ROI $B_0$. 
	Of course, one could solve the global problem
	and discard solutions that do not lie in $B_0$.
	But its complexity would not be a function of the roots in $2B_0$.
	Such local complexity behavior are provable in the univariate case 
	(e.g., \cite{becker+4:cluster:16}), and we will also
	show similar local complexity in the algorithm of this paper.
	
	\ifNeedSpace{	
	Let us briefly see the ``philosophical objections''
	of \cite{sherbrooke-patrikalakis:93}:
	The algebraic technique approach
	has two distinct phases:
	first is an algebraic (or symbolic) phase
	to reduce (e.g., using Gr\"obner bases) the function $\bff$
	to a simple form (e.g., shape lemma or triangular form);
	the second phase is numerical one that exploits the simple form to
	approximates the roots.
	The first phase restricts $\bff$ to be polynomial.
	The symbolic phase is global in nature and cannot take advantage
	of a given ROI.  On the other hand, the
	numerical phase can exploit the ROI, provided it is solved using
	subdivision-type algorithms.
	The homotopy approach needs to first solve some initial system
	$G=0$, and then follow homotopy paths from the solutions of $\bfg=\0$
	to the solutions of $\bff=\0$.  We currently have no technique 
	to predict which roots of $\bfg$ will homotopy into zeros in $B_0$.
	}
	
	Focusing on the subdivision approach, we distinguish two 
	types of subdivision: algebraic and analytic.
	In algebraic subdivision, $\bff$ is polynomial and
	one exploits representations of polynomials
	such as Bernstein form or B-splines
	\cite{mourrain-pavone:subdiv-polsys:09,
		elber-kim:solver:01,
		sherbrooke-patrikalakis:93,
		garloff-smith:bernstein-systems:01,
		garloff-smith:subdiv-algo:01}.
	Analytic subdivision
	\cite{neumaier:equations:bk-90,
		kearfott:bk,
		hentenryck+2:solving:97}
	supports a broader class of functions; this is formalized in
	\cite{yap-sagraloff-sharma:cluster:13} and includes
	all the functions obtained from composition
	of standard elementary functions or hypergeometric functions.
	Many algebraic algorithms come with complexity analysis,
	while the analytic algorithms typically lack such analysis,
	unless one views convergence analysis as a weak form
	of complexity analysis.  This lack is natural because many analytic
	algorithms are what theoretical computer
	science call ``heuristics'' with no output guarantees.
	Any guarantees would be highly\footnote{
		The issue of ``unconditional algorithms'' is
		a difficult one in analytic settings.
		Even the algorithm in this paper is conditional:
		we require the zeros of $\bff$ to be simple within $2B_0$.
		But one should certainly specify any conditions
		upfront and try to avoid conditions which are
		``algorithm-induced'' (see \cite{yap:degeneracies}).
	} conditional (cf.~ \cite{hentenryck+2:solving:97}).
    \newOld{To our knowledge, the existing subdivision algorithms, 
    both the algebraic ones and the analytic ones, suffer from a gap:
    they require an input $\vareps>0$ to serve as termination criterion
	\cite{mourrain-pavone:subdiv-polsys:09,
		elber-kim:solver:01,
		sherbrooke-patrikalakis:93,
		garloff-smith:bernstein-systems:01,
		garloff-smith:subdiv-algo:01}.
    Without this additional $\vareps$, the termination of the algorithms
    becomes unclear.
    }
    {To our knowledge, there has been no
	subdivision algorithm that solves the root isolation problem
	until the present paper.
	The subdivision algorithms
	\cite{mourrain-pavone:subdiv-polsys:09,
		elber-kim:solver:01,
		sherbrooke-patrikalakis:93,
		garloff-smith:bernstein-systems:01,
		garloff-smith:subdiv-algo:01}
	suffer from two gaps. 
	(1) Non-termination: they require an input $\vareps>0$
	to serve as termination criterion.
	(2) Non-isolation: the output box is not guaranteed to
	be \dt{isolating}, i.e., to contain a unique root. 
	So an output box could err in one of two ways:
	it may contain no roots or
	may have more than one root.  To avoid the first error,
	some root existence test is needed: so Garloff and Smith
	\cite{garloff-smith:bernstein-systems:01, garloff-smith:subdiv-algo:01}
	considered the use of Miranda test.
	To avoid the second error, Elber and Kim
	\cite{elber-kim:solver:01} introduced
	a cone test to ensure that there is at most one solution.
	The cone test generalizes the hodograph test of Sederberg and Meyers
	(1988); unfortunately this is a nontrivial test and
	details on how to compute the cones are missing.
    }
	
	\subsection{Generic Root Isolation Algorithms}
	It is useful to formulate a ``generic algorithm'' for local
	root isolation (cf.~\cite{cxy}).  We postulate 5 abstract
	modules: three box tests
	(\dt{exclusion} $C_0$, \dt{existence} $EC$, \dt{Jacobian} $JC$)
	and two box operators (\dt{subdivision} and \dt{contraction}).
	Our tests (or predicates, which we use interchangeably) 
    are best described using a notation: for any set $B\ib\RR^n$,
	$\#(B)=\#_\bff(B)$ denotes the number of roots, counted
	with multiplicity, of $\bff$ in $B$.
	These tests are abstractly defined by these implications:
		\beql{c0}
			\grouping{C_0(B)~\implies \#(B)=0,\\
			EC(B)\implies \#(B)\ge 1,\\
			JC(B)\implies \#(B)\le 1.}
		\eeql
	Unlike exact predicates, these tests are ``one-sided'' 
	(cf.~\cite{yap-sagraloff-sharma:cluster:13}) since their failure
	may have no implications for the negation of the predicate.
	For root isolation, we need both $EC(B)$ and $JC(B)$ to prove
	uniqueness.  These 3 tests can be instantiated in a variety
	of ways. The exclusion test $C_0(B)$ is instantiated differently
	depending on the type of subdivision:  
	exploiting the convex hull property of
	Bernstein coefficients (in algebraic case)
	or using interval forms of $\bff$ (in analytic case).
	For $EC$, we can use various tests coming from
	degree theory or fixed point theory
	(e.g., \cite{alefeld+3:existence:04}).  This paper
	is focused on a test based on Miranda's Theorem.
	The Jacobian test $JC$ is related to the determinant
	of the Jacobian matrix but more geometric forms
	(e.g., cone test \cite{elber-kim:solver:01}) can be formulated.
	Next consider the box operators:
	An $n$-dimensional box $B$ may be \dt{subdivided} into 
    $2^k$ subboxes in $n\choose k$ ways ($k=1\dd n$)\ifNeedSpace{, 
    giving a total of $2^n-1$ ways}.
    In practice, $k=1$ and some heuristic will choose
	one of the $n$ binary splits (see 
	\cite{garloff-smith:bernstein-systems:01} for 3 heuristics).
    \ifNeedSpace{If Bernstein form is used, then de Casteljau's algorithm
    is used to construct the Bernstein forms for the children.}
	We \dt{contract} $B$ to $B\cap N(B)$ where $N(B)$ is a box
	returned by a interval Newton-type operator.  Let us say the contraction
	``succeeds'' if the width $w(B\cap N(B))$ is less
	than $w(B)$.  But success is not guaranteed, and so this
	operator always needs to be paired with some subdivision
	operator that never fails.  It is well-known that $N(B)$
	can also provide exclusion and uniqueness tests:
		\beql{tests}
		\grouping{
		    \mbox{exclusion:} & B\cap N(B) =\es \\
		   \mbox{uniqueness:} & N(B)\ib B }.
	       \eeql
	Given the above 5 modules, we are ready
	to synthesize them into a root isolation algorithm:
	In broad outline, the algorithm maintains a queue $Q$ of
	candidate boxes.  Initially, $Q$ contains only the ROI $B_0$,
	the algorithm loops until $Q$ is empty:
	\progb{
    \lline[-3] {\sc Simple Isolate}($\bff, B_0$)
	\lline[-1] \OUTPUt: sequence of isolating boxes for roots in $B_0$
	\lline[-1] $Q\ass\set{B_0}$ 
	\lline[-1] While $Q\neq \es$
	\lline[3]  $B\ass Q.pop()$
	\lline[3]  If $C_0(B)$ continue; \Comment{discard $B$ and repeat loop}
	\lline[3]  If $EC(B)\land JC(B)$ \Comment{$B$ has a unique root}
	\lline[8] output $B$ and continue;
	\lline[3]  If $w(N(B)\cap B)<w(B)$ \Comment{if contraction succeeds}
	\lline[8] $Q$.push($B$)
	\lline[3]  else
	\lline[8] $Q$.push($subdivide(B)$)
	}

    \newOld{{\sc Simple Isolate} gives a synthetic framework of the root 
    isolating algorithms. 
    In practice, an algorithm needs not to consist of all the predicates. 
    Some of them will be sufficient.
    As mentioned above, the existing algorithms involve an input
    $\vareps$ as a criterion for termination. 
    Besides the fact that some papers lay greater emphasis on root 
    approximation than on root isolation, an important reason for 
    this phenomenon is that
    the predicates and analysis in the existing papers are not able to 
    support the termination of the algorithms without $\vareps$.
    }{
	The \dt{partial correctness} of {\sc Simple Isolate} is clear,
	i.e., if it terminates, the output is correct.
    }

    \newOld{
    For the existing algebraic subdivision algorithms, most of them have no
    existence or Jacobian test 
    \cite{sherbrooke-patrikalakis:93,mourrain-pavone:subdiv-polsys:09,
    garloff-smith:bernstein-systems:01, 
    garloff-smith:subdiv-algo:01}, others lack detailed discussion 
    on the relationship between the success of these tests and the 
    size of the boxes \cite{elber-kim:solver:01}.
    For the analytic subdivision algorithms, the interval Newton type 
    operators
    are the most favorable ones to serve as exclusion and uniqueness test.
    Extensive investigations have been performed on 
    them~\cite{neumaier:equations:bk-90,kearfott:bk}. 
    For instance, \cite[Chapter 5]{neumaier:equations:bk-90} gives detailed
    sufficient condition 
    for the strong convergence of the operators. 
    But it is still unproven that when a box is sufficiently small, the
    operators will give a definite result either to exclude the box or 
    to confirm the uniqueness of a root in it. 
    Therefore, an extra $\vareps$ is necessary 
    to ensure the termination of the algorithms. 
    But the dependence on $\vareps$ naturally results in
    two issues: the output boxes may not be isolating, 
    i.e., they may contain no root, or more than one roots.
    In this paper, we present an algorithm that makes up this gap.
    }{
    But termination is a serious issue: 
    clearly it depends on instantiations of the three tests. 
    But independent of the tests,
	non-termination can arise in two other ways:
	(1) Success of contraction ensures a reduction in the width
		$w(B)$, but this alone may not suffice for termination.
	(2) Presence of roots on the boundary of a box (e.g., $B_0$).
	We next discuss the research issues around this framework.
    }

	\subsection{How to derive effective algorithms}
	In this paper, we describe \miranda, a subdivision algorithm
	for root isolation, roughly along the above outline.  We forgo the use
	of the contraction operator as it will not figure in our analysis.  
	For simplicity, assume that all our
	boxes are hypercubes (equi-dimensional boxes); this means our
	subdivision splits each box into $2^n$ children.
	With a little more effort, our analysis can handle boxes
	with bounded aspect ratios and thus support the 
	bisection-based algorithms.
	As noted, termination depends on instantiations of our 3 tests:
	our exclusion and Jacobian tests are standard in the interval
	literature.  Our existence test, called MK test, is 
	from Moore-Kioustelides (MK) \cite{moore-kioustelidis:test:80}.
	Our algorithm is similar\footnote{
	   In \cite[Appendix]{lien-sharma-vegter-yap:arrange:14z},
	   only termination was proved (up to the abstract level)
       with no complexity analysis. 
	   We will correct an error there.
	}
	to one in the Appendix of
		\cite{lien-sharma-vegter-yap:arrange:14z}.
	In the normal manner of theoretical algorithms,
	one would proceed to ``prove that \miranda\ is correct and
	analyze its complexity''.  This will be done, but
	the way we proceed is aimed at some broader issues discussed next.

	{\bf Effectivity:} how could we convert a
	mathematically precise algorithm (like \miranda)
	into an ``effective algorithm'', i.e., certified and implementable.
	One might be surprised that there is an issue.  The non-trivially of
	this question can be illustrated from the history of isolating
	univariate roots: for about 30 years, it is known that the
	``benchmark problem''
	of isolating all the roots of an integer polynomial with $L$-bit
	coefficients and degree $n$
	has bit-complexity $\wtO(n^2L)$, a bound informally
	described as ``near-optimal''.  This is achieved by the
	algorithm of Sch\"onhage and Pan (1981-1992).
	But this algorithm has never been implemented.  What is the barrier?
	Basically, it is the formidable problem of mapping 
	algorithms in the Real RAM model \cite{bigahu}
	or BSS model \cite{bcss:bk}
	into a bit-based Turing-computable model
	-- see \cite{yap:praise:09}.
	\ignore{
	}
	\ignore{%
		   Pan (2002) [p.~705]  
	    {\small
	      ... since Sch\"onhage (1982b) already has 72 pages,
	      and Kirrinnis (1998) has 67 pages, this ruled out a
	      self-contained presentation of our root-finding 
	      algorithm.  }
		   Pan (2002) 
	    {\small
	     Our algorithms are quite involved, and their
	     implementation would require a non-trivial work,
	     incorporating numerous known implementation techniques
	     and tricks (Bini and Fiorentino, 2000; Fortune, 2001;
	     Bini and Pan, to appear). We do not touch this vast
	     domain here and just briefly comment on the precision
	     of computing.  }
	  }%

	In contrast, recent progress
	in subdivision algorithms for univariate roots
	finally succeeded in achieving comparable
	complexity bounds of $\wtO(n^2(L+n))$, and such algorithms
	were implemented shortly after!
	Thus, these subdivision algorithms were ``effective''.
	For two parallel accounts of this development, see
	\cite{sagraloff-mehlhorn:real-roots:16,
		kobel-rouillier-sagraloff:for-real:16}
	for the case of real roots, and to
	\cite{becker+3:cisolate:18,
		becker+4:cluster:16, imbach-pan-yap:ccluster:18}
	for complex roots.  What is the power conferred by subdivision?
	We suggest this:
		{\em the subdivision framework
	provides a natural way to control the numerical
	precision necessary to ensure correct operations of the algorithm.
	Moreover, the typical one-sided tests of subdivision
	avoid the ``Zero Problem'' and can be effectively implemented
	using approximations with suitable rounding modes.}

	In this paper, we capture this pathway to effectivity by
	introducing 3 Levels of (algorithmic) Abstractions:
	(A) \dt{Abstract Level},
	(I) \dt{Interval Level}, and
	(E) \dt{Effective Level}.
	We normally identify Level~(A) with the mathematical
	description of an algorithm or Real RAM algorithms.
    \newOld{At level~(I), the set extensions of the functions 
    are replaced by the interval forms (see Section~\ref{interval} for
    definitions).}{}
    \newOld{At the Effective Level, the algorithm approximate real numbers by
    BigFloat or dyadic numbers, i.e., $\ZZ[\half]$.}
    {{\em We assume our effective algorithms approximate
        real numbers by BigFloat or dyadic numbers, i.e., $\ZZ[\half]$.}}
	As illustration, consider the exclusion test $C_0(B)$
	(viewed as abstract) has correspondences in the next three levels:
	\btable{
	    (A): & $C_0(B)$ 		&$\equiv$& $0\notin \bff(B)$  \\
	    (I): & $\intbox C_0(B)$	&$\equiv$& $0\notin \intbox \bff(B)$ \\
	    (E): & $\wtbox C_0(B)$	&$\equiv$& $0\notin \wtbox \bff(B)$
	    }
	where $\bff(B)$ is the exact range of $\bff$ on $B$,
	$\intbox\bff(B)$ is the interval form of $\bff$,
    and $\wtbox\bff(B)$ the effective form\newOld{ where the endpoints 
    are dyadic numbers}{}.
	The 3 range functions here are related as
		$ \bff(B)\ib \intbox\bff(B)\ib \wtbox\bff(B). $
    \newOld{}{In general, for any abstract test $C(B)$, we derive
	its interval and effective forms to
	ensure the implications
		\beql{cb}
			\wtbox C(B)\Rightarrow \intbox C(B)\Rightarrow C(B).
		\eeql
	This means, the success of $\wtbox C(B)$ implies the success
    of $\intbox C(B)$, and hence $C(B)$.}

	An abstract algorithm $A$ is first mapped into an
	interval algorithm $\intbox A$. But the algorithm still
	involves real numbers.  So we must map $\intbox A$ to
	an effective algorithm $\wtbox A$. 
	Correctness must ultimately be shown at the Effective Level;
	the standard missing link in numerical (even ``certified'')
	algorithms is that one often stops at Abstract or Interval Levels.

    \newOld{We need to mention that the effectivity of an algorithm 
    has no implications for the efficiency of the algorithm}{}.

	{\bf Complexity:} The complexity of analytic algorithms
	is often restricted to convergence analysis.  But in this paper,
	we will provide explicit bounds on complexity as a function of the
	geometry of the roots in $2B_0$. 
	This complexity can be captured at each
	of our 3 levels, but we always begin by proving
	our theorems at the Abstract Level,
	subsequently transferred to the other levels.
	Although it is the Effective Level that really matters, it would
	be a mistake to directly attempt such an analysis at the
	Effective level: that would obscure the underlying
	mathematical ideas, incomprehensible and error prone.
	The 3-level description enforces an orderly introduction
	of new concerns appropriate to each level.  Like structured
	programming, the design of effective algorithms needs
	some structure.  Currently, outside of the subdivision framework,
	it is hard to see a similar path way to effectivity.

	\subsection{Literature Survey}
	There is considerable literature associated with each of our three
	tests: the exclusion test comes down to bounding
	range of functions, a central topic in Interval Analysis
	\cite{ratschek-rokne:range:bk}.  The Jacobian test is
	connected to the question of local injectivity of functions, the
	Bieberbach conjecture (or de Branges Theorem), Jacobian Conjecture, and
	theory of univalent functions.  In our limited space, we focus on the
	``star'' of our 3 tests, i.e., the existence test. 
	It is the most sophisticated
	of the 3 tests in the sense that some nontrivial global/topological
	principle is always involved in existence proofs. 
	In our case, the underlying principle is the fixed point theorem
	of Brouwer, in the form of Miranda's Theorem (1940), and 
	intimately related to degree theory. 
	
	We compare two box tests $C$ and $C'$ in terms of
	their relative \dt{efficacy}:
	say $C$ is \dt{as efficacious as} $C'$,
	written $C\succeq C'$, if for all $B$,
	$C'(B)$ succeeds implies that $C(B)$ succeeds.
	The relative efficacy of several existence tests have been studied
		\cite{goldsztejn2007comparison,
			frommer2004comparison,
			frommer-lang:miranda-type-tests:05,
			alefeld+3:existence:04}.
	Goldsztejn considers four common existence tests,
	and argues that ``in practice'' there is an efficacy hierarchy
		\beql{in}
		(IN) \succeq  (HS) \succeq (FLS) \succeq (K)
		\eeql
	where
	(K) refers to Krawcyzk,
	(HS) to Hansen-Sengupta,
	(FLS) to Frommer-Lang-Schnurr,
	and (IN) to Interval-Newton.
	Note that (K), (HS) and (IN) are all based on
	Newton-type operators (see \refeq{tests}).
	Our Moore-Kioustelidis (MK) test is essentially (FLS).
	We say ``essentially'' because the details of defining
	the tests may vary to render the comparisons invalid. 
	In our MK tests, we evaluate $\bff$ on each box face using the
	Mean Value Form expansion at the center of the face. 
	But the above analysis assumes an expansion is at the center
	of the box, which is less accurate.  But we may
	also compare these tests in terms of their complexity
	(measured by the worst case number of
	arithmetic operations, or number of function evaluations);
	a complexity-efficacy tradeoff may be expected.
    \newOld{}{Such complexity comparisons do not account for adaptive costs:
	Newton-type existence tests have non-adaptive costs
	while the Miranda-type tests are adaptive (we are testing
	$n$ pairs of faces, and can break off as soon as one
    pair fails the test.}
    Finally, evaluating these tests 
	in isolation does not tell us how they might perform
	in the context of an algorithm.  It is therefore
	premature to decide on the best existence test.
	
	\ignore{
	NOTE: The algorithm introduced by this paper was first
	sketched in an Appendix in
	\cite{lien-sharma-vegter-yap:arrange:14z}.
	Lemma XXX in
	\cite{lien-sharma-vegter-yap:arrange:14z}
	was deficient; moreover, only termination
	at the A-level was proved.  In this paper,
	we prove not only termination but give
	quantitative bounds for the E-level of abstract.
	}

	\ignore{%
	In some papers, subdivision method is used along with reduction
	methods. If the reduction ratio is not as large as expected, 
	then a subdivision step is performed. In
	\cite{sherbrooke-patrikalakis:93,mourrain-pavone:subdiv-polsys:09},
	the polynomial systems are 
	represented with Bernstein basis, and control points are explored to
	design reduction step. In \cite{kearfott:bk,hansen:bk},
	algorithms based on the interval Newton methods are described.
	
	One problem is that they mostly has an $\eps$-termination
	criteria based on an input parameter $\eps$.
	
	The predicates in this paper are of three sorts:
	The key predicate is a well-known test 
	from Moore and Kioustelides
	\cite{neumaier:equations:bk-90,moore-kioustelidis:test:80},
	which is in turn based on the Poincare-Miranda Theorem.
	Although this is well-known in the interval community, their
	integration into a complete algorithm for root isolation seems
	to be missing in the literature.
	Other ingredients are an exclusion test and the Jacobian test
	(see \cite{lien-sharma-vegter-yap:arrange:14z,
	bennett-papadopoulou-yap:minimization:16}
	where similar combinations of these tools have been used).
	
	Following \cite{yap-sagraloff-sharma:cluster:13,becker+4:cluster:16},
	we generalize the problem of root isolation to natural root clustering.
	This allows our algorithm to work with polynomials whose
	coefficients are only number oracles, and even analytic functions.
	}%

\subsection{Overview}
	In section 2, we introduce some basic concepts of interval arithmetic
	and establish notations. 
	Section 3 introduces the key existence test based on Miranda's theorem.
	Section 4 proves conditions that ensure the success of these existence
	test.
	Section 5 introduces two Jacobian tests.
	Section 6 describes our main algorithm.
	Section 7 is the complexity analysis of our algorithm.
	We conclude in Section 8.
	All proofs are relegated to the Appendix.

\section{Interval Forms}\label{interval}
	We first establish notations for
	standard concepts of interval arithmetic.
	Bold fonts indicate vector variables: e.g.,
	$\bff=(f_1\dd f_n)$ or $\bfx=(x_1\dd x_n)$.

	Let $\intbox \RR$ denote the set of compact intervals in $\RR$.
	Extend this to $\intbox\RR^n$ for the set of compact $n$-boxes.
	In the remaining paper, we assume that all $n$-boxes are hypercubes
	(i.e., the width in each dimension is the same).
	For any box $B\in\intbox\RR^n$, let $\bfm_B= \bfm(B)$ denote its center
	and $w_B=w(B)$ be the width of any dimension.
	Besides boxes, we will also use ball geometry:
	let $\Delta=\Delta(\bfa,r)\ib\RR^n$ denote the closed
	ball centered at $\bfa\in\RR^n$ of radius $r>0$.
	If $r\le 0$, $\Delta(\bfa,r)$ is just the point $\bfa$.
    \newOld{For any $k>0$, let $kB$ denote
    the box centered at $\bfm(B)$ of width $k\cdot w(B)$, called
    the \dt{$k$-dilation} of $B$.
    The $k$-dilation $k\Delta$ of $\Delta$ is defined likewisely.}
    {For any positive $k>0$, let $k\Delta$ and $kB$
	denote the dilation of the ball $\Delta$ and box $B$ relative
	to their centers.  }
    
	Let $A,B\ib\RR^n$ be two sets.  We will 
	quantify their ``distance apart'' in two ways:
	their usual Hausdorff distance is denoted $q(A,B)$
	and their \dt{separation},
    	 $\inf\set{\|\bfa-\bfb\|:
        \bfa\in A, \bfb\in B}$ is denoted as $\sep(A,B)$.
	Note that $q$ is a metric on closed subsets of $\RR^n$
	but $\sep(A,B)$ is no metric.

	Consider two kinds of extensions of a
	function $f:\RR^n\to \RR$.  First, the
	\dt{set extension} of $f$ refers to the function
	(still denoted by $f$) that maps
	$S\ib \RR^n$ to $f(S)\as \set{f(\bfx): \bfx\in S}$.
	The second kind of extension is not unique:
	an \dt{interval form} of $f$ is any function 
		$\intbox f: \intbox \RR^n\to \intbox \RR$,
	satisfying two properties:
	(i) (inclusion) $f(B)\ib \intbox f(B)$;
    \newOld{(ii) (convergence) 
    if $B_1\supseteq B_2\supseteq \cdots \supseteq B_i
    \supseteq \cdots$
    with $\bfp=\lim_{i=0}^\infty B_i$
    then  $f(B_1)\supseteq f(B_2)\supseteq \cdots \supseteq f(B_i)\supseteq
    \cdots$ and
    $f(\bfp)=\lim_{i=0}^\infty \intbox f(B_i)$.}
    {
    if $\bfp=\lim_{i=0}^\infty B_i$
	then $f(\bfp)=\lim_{i=0}^\infty \intbox f(B_i)$.  }
	For short, we call $\intbox f$ a \dt{box form} of $f$.
	If $\bff=(f_1\dd f_n):\RR^n\to\RR^n$, we have corresponding
	set extension $\bff(S)$ 
	and interval forms $\intbox \bff:\intbox \RR^n\to\intbox\RR^n$.

    \newOld{}
    { For any set $S\ib\RR^n$, 
	let $\zero_\bff(S)$ denote the multiset of zeros of $\bff$ in $S$.
	We assume that $\bff$ is analytic and its zeros are
	counted with the proper multiplicity.  
	Then $\#_\bff(S)$ is the size of the multiset $\zero_\bff(S)$.
	We may write $\zero(S)$ and $\#(S)$ when $\bff$ is understood.
    }
	
	The notation ``$\intbox f$'' is generic;
	we use subscripts to indicate specific box forms.
	Thus, the \dt{mean value form} of $f$ is 
		$$\boxm f(B) =
			f(\bfm(B))+\intbox\nabla f(B)^T\cdot (B-\bfm(B))$$
	where $\nabla f$ is the gradient of $f$ (viewed as a column vector)
	and $\nabla f(B)^T$ is the transpose.  The box $B-\bfm(B)$ is
	now centered at the origin, i.e., $\bfm(B-\bfm(B))=\0$.
	The appearance of the generic
	``$\intbox\nabla f(B)$'' in the definition
	of $\boxm f$ means that $\boxm f$ is still not
	fully specified.  In our complexity analysis, we assume that
	for any box form, if not fully
	specified, will have at least linear convergence.
    \newOld{}{{\em In this paper, all the box forms used in our
    predicates will be mean value forms.}}
        Next, we intend
	to convert the interval form $\boxm$ to some effective version
	$\wtboxm$.  
    \newOld{}{One reason that this is necessary may be seen in
	the fact that $\boxm$ assumes an exact value $f(\bfm(B))$.  Even
	if $\bfm(B)$ is a dyadic number, we may need to approximate
    $f(\bfm(B))$ (e.g., $f(x)=\sin(x)$).}


\section{Miranda and MK Tests}
	In the rest of this paper, we fix 
		\beql{f}
		    \bff \as (f_1\dd f_n): \RR^n \to \RR^n
		   \eeql
	to be a $C^2$-function (twice continuously differentiable),
	and $\bff$ and its partial derivatives have interval forms.
	We further postulate that $\bff$ has only finitely many
	simple zeros in $2B_0$ where $B_0$ is the bounded
    region of interest.
	A zero $\bfalpha$ of $\bff$ is simple if
	the Jacobian matrix $J_\bff(\bfalpha)$ is non-singular.
    \newOld{For any set $S\ib\RR^n$, 
	let $\zero_\bff(S)$ denote the multiset of zeros of $\bff$ in $S$.
	We assume that $\bff$ is analytic and its zeros are
	counted with the proper multiplicity.  
	Then $\#_\bff(S)$ is the size of the multiset $\zero_\bff(S)$.
    We may write $\zero(S)$ and $\#(S)$ when $\bff$ is understood.}{}
    The \dt{magnitude} of
	any bounded set $S\ib\RR$ is defined as 
    $|S| \as\sup\set{|x|: x\in S}$.

	We consider a classical test from Miranda (1940)
	to confirm that a box $B\in\intbox\RR^n$ contains a
	zero of $\bff$.
	If the box $B$ is written as $B=\prod_{i=1}^n I_i$ with
	$I_i = [a_i^-,a_i^+]$, then it has two $i$-th \dt{faces}, namely
		$$B^-_i \as
			I_1\times \cdots \times I_{i-1}\times
            		\set{a^-_i}\times I_{i+1}\times\cdots \times I_n$$
	and $B^+_i$, defined similarly.  Write $B^{\pm}_i$
	to mean either $B^-_i$ or $B^+_i$.  
	Consider the following box predicate called\footnote{
		We call it ``simple'' as we ignore some common
		generalizations
		that allow an interchange of ``$<0$'' with ``$>0$'',
		or replace $\bff$ by
			$\sigma(F)=(f_{\sigma(1)}\dd f_{\sigma(n)})$ for any
		arbitrary permutation $\sigma$ of the indices.
	}
	the \dt{simple Miranda Test}:
	       \beql{smt}
	       \mtf(B)\equiv\quad
			\doublebigwedge\nolimits_{i=1}^n
				(f_i(B_i^+) > 0) \land (f_i(B_i^-) < 0)
		\eeql
	    where $\bff$ is given in \refeq{f}.
	   The following result is classic:
	\bprol{simpleMiranda}{[Miranda (1940)]} \ \\
	    \hspace*{15mm} If $\mtf(B)$ holds then $\#_\bff(B)\ge 1$.
	\eprol
    \newOld{}{For a box $B$ and $k>0$, let $kB$ denote
    the box centered at $\bfm(B)$ of width $k\cdot w(B)$, called
    the \dt{$k$-dilation} of $B$.}
    Next, we introduce the \dt{MK Test} test $\mk(B)=\mkf(B)$ 
    that amounts an application of the simple Miranda test to the
	box $2B$, using a preconditioned form of $\bff$:
    \progb{
	\lline[-3] {\sc Abstract MK Test}
	\lline[0] \INPUt: $\bff$ and box $B$
	\lline[0] \OUTPUt: \true\ iff $\mkf(B)$ succeeds
	\lline[5] 
        \newOld{1. Take a point $\bfm\in B$ with
            $J\inv(\bfm)$ well-defined.}{}
	\lline[5] 2. Construct a ``preconditioned version'' $\bfg$:
	\lline[13]
            $\bfg\ass J\inv(\bfm) \bff = (g_1(\bfx)\dd g_n(\bfx))$
	\lline[5] 3. Apply the Simple Miranda Test to $\bfg$ over $2B$:
	\lline[13] For $i\ass 1\dd n$:
	\lline[18] 	If $g_i(2B^+_i)\le 0$ or $g_i(2B^-_i)\ge 0$,
			\hspace*{7mm} (*)
	\lline[23] return \false
	\lline[5] 4. Return \true.
	}
    
	The notation ``$2B^{\pm}_i$'' in (*) 
	refers to faces of the box $2B$, not the $2$-dilation of the
	faces of $B$.  Here ``MK'' refers to Moore and Kiousteliades
	\cite{moore-kioustelidis:test:80};
	the preconditioning idea first appearing in
	\cite{kioustelidis:miranda:78}.
	The MK Test was first introduced in
	\cite{lien-sharma-vegter-yap:arrange:14z}.
    \newOld{Notice that in \mk($B$), the Miranda test is performed on
    $2B$ instead of $B$. It is intended to address the difficult case 
    where the root is close to the boundary of a box.}{}

	Note that \mk($B$) is mathematically exact and
	generally not implementable (even if it were possible,
	we may still prefer approximations).
	We first define its interval form, denoted
	$\intmk(B)$: simply by replacing 
    $g_i(2B^{\pm}_i)$ in line (*) by interval forms 
    $\intbox g_i(2B^{\pm}_i)$.
	Finally, we must define the effective form
	$\wtmk(B)$ (Section 8).  
    \newOld{}{The key property is the relation
	(cf.~\refeq{cb}):
        $$ \wtmk(B) \Rightarrow \intmk(B) \Rightarrow \mk(B).$$}

\section{On Sure Success of MK Test}
	The success of the MK test implies the existence of roots.
	In this section, we prove some (quantitative) converses.

	We need preliminary facts about mean value forms.
	Given $x,y\in \RR$, the notation $x\pm y$ denotes a number of the form
	$x + \theta y$, where $0\le |\theta| \le 1$; thus
	``$\pm$" hides the implicit $\theta$ in the definition.
	This notation is not symmetric: $x\pm y$ and $y\pm x$
	are generally different.
	This notation extends to matrices:
	let $A=(a_{ij})_{i,j=1}^n$ and $B=(b_{ij})_{i,j=1}^n$ be two matrices,
	then $A\pm B \as (a_{ij} \pm b_{ij})_{i,j=1}^n$.
    \newOld{}{Similarly, for a scalar $\lambda$, we have
    $A\pm\lambda\as (a_{ij}\pm \lambda)_{i,j=1}^n$.}
	Also, let $|\bfx|$ denote the vector $(|x_1|\dd |x_n|)$ where
	$\bfx=(x_1\dd x_n)$.
	For $\bfx,\bfy\in \RR^n$, we write $[\bfx,\bfy]$
	to denote the line segment connecting $\bfx$ and $\bfy$.
	We write $\|\bfx\|$ and $\|A\|$
	for the infinity norms of vector $\bfx$ and matrix $A$.
    For a bounded convex set $C\ib\RR^n$,
	define the matrix $K(C)$ with entries $(K(C)_{ij})_{i,j=1}^n$
	where
		\beql{kij}
			K(C)_{ij} \as \sum\nolimits_{k=1}^n \Big|
			\frac{\partial^2 f_i}{
		    	\partial x_j\partial x_k} (C)\Big|.
		\eeql
	Below, $C$ may be a disc $\Delta$ or a line $[\bfx,\bfy]$.
    Denote by $J_{\bff}(\bfx)$ the Jacobian matrix of $\bff$ at $\bfx$.
    We write $J_{\bff}(\bfx)$ as $J(\bfx)$ when $\bff$ is understood.
	The following is a
	simple application of the Mean Value Theorem (MVT):
	
	\blemT{MVT}{mvt}
	Given two points $x,y\in \RR^n$, we have:
    \begin{flalign}
        &\text{(a)} \quad\quad \quad\quad 
        J(\bfx) = J(\bfy) \pm \K{[\bfx,\bfy]}\|\bfx-\bfy\|,& \nonumber \\
        &\text{(b)}\ \bff(\bfx) - \bff(\bfy) = (J(\bfy) \pm \K{[\bfx,\bfy]}
            \|\bfx-\bfy\|)\cdot(\bfx-\bfy).& \nonumber
    \end{flalign}
	\elemT
	
\subsection{Sure Success of abstract MK Test}
    In this and the next subsection,
    we consider boxes that contain a root
    $\bfalpha$ of $\bff$.  We prove conditions that ensures
    the success of the MK Test.  We first prove this
    for the abstract test $\mk(B)$.  The next section
    extends this result to the interval test $\intmk(B)$.

    The key definition here
    is a bound $\lambda_1(\bfalpha)$ which depends on
    $\bfalpha$ and $\bff$.  We prove that if $w(B)\le \lambda_1(\bfalpha)$,
     then the abstract MK test will succeed on $B$.
     By a \dt{critical point} we mean $\bfa\in\RR^n$ where
     the determinant of $J(\bfa)$ is zero.  
     By definition, a
     root $\bfalpha$ of $\bff$ is simple if $\bfalpha$ is not a critical point.

     Suppose $S_1$ and $S_2$ are two \newOld{bounded}{} sets in $\RR^n$.
     Define
       \newOld{$$\|J\inv(S_1)\|\as \sup\nolimits_{\bfx\in S_1}
            \|J\inv(\bfx)\|\quad \text{and} $$}
       {    $$\|J\inv(S_1)\|\as \max\nolimits_{\bfx\in S_1}
            \|J\inv(\bfx)\|\quad \text{and} $$}
       \newOld{$$\|J\inv(S_1)\cdot K(S_2)\|
	    	\as \sup\nolimits_{\bfx\in S_1, \bfy\in S_2}
            \|J\inv(\bfx)\cdot K(\bfy)\|. $$}
       {    $$\|J\inv(S_1)\cdot K(S_2)\|
	    	\as \max\nolimits_{\bfx\in S_1, \bfy\in S_2}
            \|J\inv(\bfx)\cdot K(\bfy)\|. $$}
    We see that both $\|J\inv(S_1)\|$ and $\|J\inv(S_1)\cdot K(S_2)\|$
    are finite if $S_1$ does not contain a critical point of $\bff$.
    \ignore{ 
	Let $\Delta=\Delta(\bfa,r)\ib\RR^n$ denote the closed
	disc centered at $\bfa\in\RR^n$ of radius $r$.
	For $r\le 0$, define $\Delta(\bfa,r)$ to be the point $\bfa$.
	For any real $k$, let $k\Delta$ denote $\Delta(\bfa,kr)$.
	}%
    Consider the following function
	\beql{radius}
       s(r)\as r - \efrac{27n\jkd{\bfalpha,2\sqrt{n}r}{\bfalpha,2\sqrt{n} r}}.
	\eeql
	We then define $\lambda_1(\bfalpha)$ to be the smallest $r$
	such that $s(r)=0$, i.e.,
	$\lambda_1(\bfalpha)\as \argmin_r \set{ s(r)=0 }$.
	
	\blem
	For any simple root $\bfalpha$ of $\bff$, $\lambda_1(\bfalpha)$
	is well-defined.
	\elem
    
	From now on, let $\discalpha$ denote the disc
    \beql{delta1}
    \discalpha \as \Delta(\bfalpha,2\sqrt{n}\lambda_1(\bfalpha)).
    \eeql
	The following lemma corrects an gap in the appendix of
	\cite{lien-sharma-vegter-yap:arrange:14z}.

	\bleml{theoremMK}
    \newOld{
        Let $B$ be a box containing a simple root $\bfalpha$ of 
        $\bff$ and $\bfm \in B$ with $J\inv(\bfm)$ well-defined. 
        If $w_B\le \lambda_1(\bfalpha)$, the preconditioned system
	    $\bfg_B \as J\inv (\bfm)\bff = (g_1\dd g_n)$ satisfies that
    }{  Let box $B$ contain a simple root $\bfalpha$ of $\bff$.
	    \\ If $w_B\le \lambda_1(\bfalpha)$,
	    the preconditioned system
	    $\bfg_B \as J\inv (\bfm(B))\bff = (g_1\dd g_n)$ is
	    well-defined, and 
    }
        for all $i=1\dd n$,
	    	$$g_i(2B^+_i) \ge \frac{w_B}{4},\qquad
	    	g_i(2B_i^-) \le -\frac{w_B}{4}.$$
	\eleml

\subsection{Sure Success of Interval MK Test }
    We now extend the previous subsection on
    the abstract MK Test $\mk(B)$ to the interval version $\intmk(B)$. 
    Again, assume $B$ is a box containing exactly one
    root $\bfalpha$ of $\bff$.
    We will give $\lambda_2(\bfalpha)$ which is analogous to 
    $\lambda_1(\bfalpha)$ and prove that if 
    $w_B\le \lambda_2(\bfalpha)$, then
    $\intmk(B)$ will succeed.

    To prove the existence of such a $\lambda_2(\bfalpha)$ as
    mentioned above, we need to make some assumptions on the
    property of the box functions.
    As in \cite{moore:bk-95}, a box function
     $\intbox f$ is called \dt{Lipschitz} in a region $S\ib\RR^n$
    if there exists a  constant $L$ such that
    \beql{lipschitz}
    w(\intbox f(B))\le L\cdot w(B), \quad \forall B\ib S.
    \eeql
    We call any such $L$  a \dt{Lipschitz constant} of $\intbox f$ on $S$.
	For our theorem, we need to know the specific box function
	in order to derive a Lipschitz constant.  Consider the
    	mean value form $\boxm f$ on a region $S\ib \RR^n$.
    \blem
    Let $f$ be a continuously differentiable function defined on a 
    convex region $S\ib \RR^n$.  
    Then 
    $\sum_{k=1}^n \Big|\intbox \frac{\partial f}{\partial x_j}(S)\Big|$
    is a Lipschitz constant for $\boxm f$ on $S$.
    \elem
    
    Consider the sign tests of $\intmk(B)$:
        $$\boxm g_i(2B^+_i)>0 \quad \text{and} \quad
		\boxm g_i(2B^-_i)<0$$
    where $g_i$ is the $i$-th component of the system $J(\bfm))\inv \bff$.
    We consider the mean value form
        $\boxm g_i(2B^+_i) = g_i(m(2B^+_i))+\intbox
        \nabla g_i(2B^+_i)\cdot (m(2B^+_i)-2B^+_i)$
    and assume that the components of $\intbox\nabla g_i(2B^+_i)$
    are evaluated via the linear combination of
    $\intbox \frac{\partial f_j(2B^+_i)}{\partial x_k}$
    for $j,k=1\dd n$.

    We now prove that if $B$ is small enough,
    $\intmk(B)$ will succeed.  Recalling the Hausdorff distance $q(I,J)$
    on intervals, we have this bound from \cite{neumaier:equations:bk-90}.
	\bprol{quadConverg}
	   Let $f:D\subset \RR^n \rightarrow \RR$ be a continuously
	   differentiable function.
	   Then
	      \beq
	        q(\boxm f(B), f(B)) \le 2w_B \sum\nolimits_{i=1}^n
			w(\intbox\frac{\partial f(B)}{\partial x_i}).
	      \eeq
	\eprol

	For the next theorem, define
		\beql{whlambda1}
			\wh{\lambda}_1(\bfalpha)\as
			\efrac{64n^2 L\cdot \|J\inv(\discalpha)\|}.
		\eeql
	where $L=L_\bfalpha$ is a Lipschitz constant for
        	$\intbox \frac{\partial f_j}{\partial x_k}$
		on $\discalpha$  (for all $j,k=1\dd n$).
	\bthml{MK}
	    Let $B$ be a box containing a simple root $\bfalpha$
        of width $w_B\le \lambda_{1}(\bfalpha)$ and 
        $\bfm\in B$ with $J\inv(\bfm)$ well-defined.
        \benum[(a)]
        \item 
	    If $w(\intbox\frac{\partial g_i(2B_i^+)}{\partial x_j})
			\le \frac{1}{32n}$ for each $j=1\dd n$, 
        then $\intmk(B)$ succeeds with $\bfg_B := J\inv (\bfm)\bff$.
        \item
        	If $w_B\le \lambda_2(\bfalpha)$ with
        	$\lambda_2(\bfalpha)\as\min\set{\lambda_1(\bfalpha),
        		\wh{\lambda}_1(\bfalpha)}$,
        	then $\intmk(B)$ succeeds.
        \eenum
	\ethml
\section{Two Jacobian Conditions}
	We define the \dt{Jacobian test} as follows:
        	\beql{jc}
		\jc(B)\equiv	0 \notin \det(J_\bff(3B)). 
		\eeql
	The order of operations in $\det (J_\bff(3B))$
	should be clearly understood: first we compute the
	\dt{interval Jacobian matrix} $J_\bff(3B)$, i.e., 
	entries in this matrix are the intervals $\partial_{x_j}f_i(3B)$.
	Then we compute the determinant of the interval matrix.
	Also note that we use $3B$ instead of $B$.
	The following is well-known in interval computation
	(see \cite[Corollary to Theorem 12.1]{aberth:precise:bk}):

	\bprol{jacobian}{[Jacobian test]} \ \\
	    \hspace*{15mm}
	    If $\jc(B)$ holds then $\#_\bff(3B)\le 1$.
	\eprol
	\ignore{
	In fact, without involving the MK test,
	$\jc(B)$ alone can be used to verify uniqueness provided that 
	the determinant is evaluated in a specific way. 
    	It is easy to see that $\jc(B)$ alone does not ensure a unique root.
	}
	We next introduce the following \dt{strict Jacobian test}:
		\beql{jsw}
		\jcs(B)\equiv 0 \notin (\det J_\bff)(3B)
		\eeql
	where $(\det J_\bff)(\bfx)$ denotes the expression obtained
	by computing the determinant of the Jacobian matrix $J_\bff(\bfx)$
    with functional entries $\frac{\partial{f_i}}{\partial x_j}(\bfx)$.  
    Finally,
	we evaluate $(\det J_\bff)(\bfx)$ on $3B$.
	Note that $\jc(B)\Rightarrow \jcs(B)$ and so the strict test
	is more efficacious.  Unfortunately, it is known that
	$\jcs(B)$ does not imply $\#_\bff(3B)\le 1$.
	Nevertheless, we now show that it can serve
	as a uniqueness test in conjunction with the MK test:
	    \bthml{weakjc}\ \\
	    \hspace*{5mm}
	        If both $\jcs(B)$ and
		$\mk(\frac{3}{2}B)$ succeed then $\#_\bff(3B)=1$.
	    \ethml
	It follows that we could use $\jcs(B)\land \mk(B)$
	in our \miranda\ algorithm in the introduction.

\section{The \miranda\ Algorithm}
    Our main algorithm for root isolation is given
    in Figure \ref{fig:miranda}.  We use $\mk(B)$ and $\jc(B)$
    (respectively) for its existence and Jacobian tests.
    It remains to specify the exclusion test $\cz(B)$: 
		\beql{cz}
		\cz(B)\equiv (\exists i= 1\dd n)[ 0\notin f_i(B)]
		\eeql
	\begin{figure}
	\myprogb{
	    ~\\
	    {\sc Abstract \miranda($\bff, B_0$)}\\
	    {\sc Output:} Queue $P$ of non-overlapping 
	    isolating boxes of $\bff$ s.t.\\
	     \>\>\>\>   $\calZ_\bff(B_0)\subseteq$ 
            $\bigcup_{B\in P}\zero_\bff(B)\ib \zero_\bff(2B_0)$ \\
	    1. Initialize output queue $P\ass \emptyset$
		and priority queue $Q \ass \set{B_0}$. \\
	    2. While $Q\neq\emptyset$ do:\\
	    3.\> Remove a biggest box $B$ from $Q$. \\
	    4.\> If $\cz(B)$ succeeds, continue; \\
	    5.\> If $\jc(B)$ succeeds then \\
	    6.\>\> Initialize new queue $Q'\ass \set{B}$. \\
	    7.\>\> While $Q'\neq\emptyset$ do:\\
	    8.\>\>\> $B'\ass Q'.pop()$. \\
	    9.\>\>\> If $\cz(B')$ fails then \\
	    10.\>\>\>\> If $\mk(B')$ succeeds then\\
	    11.\>\>\>\>\> $P.add(2B')$. \\
	    12.\>\>\>\>\> Discard from $Q$ the boxes contained in $3B$. \\
	    13.\>\>\>\>\> Break.\\
	    14.\>\>\> $Q$.push($subdivide(B)$). \\
	    15.\> Else \\
	    16.\>\> $Q$.push($subdivide(B)$).
	}
	    \caption{Root Isolation Algorithm}
	    \label{fig:miranda}
	\end{figure}
    	The algorithm in Figure \ref{fig:miranda} is abstract.
	To introduce the interval version \intbox\miranda,
	just replace the abstract tests by their 
	interval analogues: $\intmk(B), \intcz(B)$ and $\intjc(B)$.
	It amounts to replacing the set theoretic function in the abstract
	definition by their interval analogues:
	    \bitem 
	    \item $\intcz(B)$: $\exists i= 1\dd n$ such that
		    $0\notin \intbox f_i(B)$;
	    \item 
		$\intjc(B)$: $0 \notin \intbox \det (J(3B))$;
	    \item
		In the definition of $\mk(B)$ (Section 3), replace
		each $g_i(2B^{\pm}_i)$ by
		$\intbox g_i(2B^{\pm}_i)$.
	    \eitem
    \newOld{}{Note that all these box forms are really
    mean value forms $\boxm$.}
	For the effective version, we
	use the tests $\wtmk(B), \wtcz(B)$ and $\wtjc(B)$,
	which will be discussed in Section 8.
	
	Termination of each version of
	\miranda\ follows from the complexity analysis below.
    \newOld{}{Even if there are roots on the
	boundary of $B_0$, we will terminate, although the isolated 
    root might lie in $2B_0\setminus B_0$.}
	We first show the partial correctness:
	\bthmT{Partial Correctness}{partial}
    \newOld{
        If \miranda\ halts, the output queue $P$ is correct.
    }{
	    \ \\1. If \miranda\ halts, the output queue $P$ is correct.
	    \ \\2. The same holds for \intbox\miranda\ and $\wtbox$\miranda.
    }
	\ethmT
	

\section{Complexity Upper Bounds}\label{complexity}
	In this section, we derive a lower bound $\lambda>0$
	on the size of boxes produced by \miranda.  
	That is, any box $B$ with width $w(B)\le \lambda$ will 
	either be output or rejected.  
	This implies that the subdivision tree is no deeper than
	$\log_2( w(B_0)/\lambda)$, yielding an
	upper bound on computational complexity.
	This bound $\lambda$ will be expressed in terms of quantities
	determined by the zeros in $2B_0$.
	We first prove this for the abstract \miranda,
	then extend the results to $\intbox\miranda$
	and $\wtbox\miranda$.
    From the algorithm, we see that a box $B$ is output if
    $\neg \cz(B)\wedge \jc(B)\wedge \mk(B)$
    holds in line 10;
    it is rejected if one of the 2 following cases is true:
     (1) $\cz(B)$ holds or
      (2) it is contained in $3B'$ where $B'$ is an output box, 
      as indicated in line 12.
    The boxes that contain a root
	of $\bff$ will be finally verified by the former predicate and
	the boxes that contain no root of $\bff$ will eventually be
    rejected in one of the 2 cases.

    To prove the existence of such a $\lambda$, we need to look into
    the tests $\cz(B)$, $\jc(B)$ and $\mk(B)$.
    We will give bounds  $\lambda_{\jc}$, $\lambda_{\mk}$ and
    $\lambda_{\cz}$ for the $3$ tests respectively and show that
    for any box $B$ produced in the algorithm \\
    (1) if $\#(B)>0$, it will pass $\mk(B)$ when 
     	$w_B\le \lambda_{\mk}$,\\
    (2) if $\#(B)\le 1$, it will pass $\jc(B)$ when 
    	$w_B \le \lambda_{\jc}$;\\
    \newOld{
    (3) if $\#(B) = 0$, there are 2 cases:
        (a) if $B$ keeps a certain distance from the roots,
    	it passes $\cz(B)$ when $w_B\le \lambda_{\cz}$;
        (b) if $B$ is close enough to the roots, it will be
        rejected by line 12 of the algorithm when 
        $w_B \le 
        \efrac{2}\min\{\lambda_{\jc},\lambda_{\mk},\lambda_{\cz}\}$.}
    {(3) if $\#(B) = 0$ and $B$ keeps a certain distance from the roots,
    	it will pass $\cz(B)$ when $w_B\le \lambda_{\cz}$.
    }
        \\
    We have essentially proved item (1) in the Section 4.
    More precisely, for each root $\bfalpha$, we had defined a constant
    $\lambda_2(\bfalpha)$.  We now set
	\beql{lmk}
	\lambda_\mk \as \min_{\bfalpha\in \zero(2B_0)} \lambda_2(\bfalpha).
	\eeql

    \subsection{Sure Success for $\cz(B)$ and $\jc(B)$}
	We study conditions to ensure the
	success of the tests $\jc$ and $\cz$.
	We will introduce constants $\lambda_\jc, \lambda_\cz$ in
	analogy to \refeq{lmk}.
	
	First consider $\jc(B)$.
	\ignore{%
	Let $\calS_J$ be the solution set of $\det(J(\bfx))=0$
	for $\bfx\in B_0$, and define
	$d_1\as \sep(\zero(B_0),\calS_J)$.
	Note that $d_1>0$ since
	the roots of $\bff$ in $B_0$ are simple.
	For any box $B$ containing a root, we
	claim that $w_B < \frac{d_1}{2\sqrt n}$ implies
	$3B\cap\calS_J = \emptyset$.
	[In proof: $q(\bfalpha,\del (3B))< 2w_B\sqrt{n}$.
	Hence $3B\ib \Delta(\bfalpha,2w_B\sqrt{n})\ib \Delta(\bfalpha,d_1)$.]
	Thus $0\notin\det(J(B))$, and the condition $\jc(B)$ holds.
    {\coblue{There is an error if we use the strong Jacobian condition.
    How to correct it? Use the same argument as in the interval version?}}
    }
    Let box $B$ contain a simple root $\bfalpha$.
    By Mean Value Theorem, 
    $w(\frac{\partial f_i}{\partial x_j}(3B))\le3w_B\cdot K(3B)_{ij}$
    (see \refeq{kij} for definition).
    Since 
    $\frac{\partial f_i}{\partial x_j}(\bfalpha)\in
    \frac{\partial f_i}{\partial x_j}(3B)$,
    it holds 
    $\frac{\partial f_i}{\partial x_j}(3B)\ib
        [\frac{\partial f_i}{\partial x_j}(\bfalpha)-3w_B\cdot K(3B)_{ij},
        \frac{\partial f_i}{\partial x_j}(\bfalpha)+3w_B\cdot K(3B)_{ij}]$
   ($\forall i,j=1\dd n$).
   Denoting $U\as \max_{1\le i,j\le n}
   	|\frac{\partial f_i}{\partial x_j}(\bfalpha)|$
   and 
   $V \as$ $\max_{1\le i,j\le n }\cdot K(3B)_{ij}$,
   we get 
    $|\frac{\partial f_i}{\partial x_j}(3B)| \le U+3Vw_B$
    and 
    $w(\frac{\partial f_i}{\partial x_j}(3B)) \le3Vw_B$.
   By applying the rules 
   $w(I_1+I_2)=w(I_1)+w(I_2)$ and
   $w(I_1\cdot I_2)\le w(I_1)\cdot |I_2|+w(I_2)\cdot |I_1|$
   where $I_1$, $I_2$ are intervals,
   we may verify by induction that
   $w(\prod_{i=1}^n(\frac{\partial f_i}{\partial x_{\sigma_i}}(3B))\le 
   3nV(U+3w_BV)^{n-1}w_B$ 
   for any permutation $\sigma$.
   Hence, it follows 
   $w(\det(J_\bff(3B)))\le 3n\cdot n!\cdot V(U+3Vw_B)^{n-1}w_B$.

   Set $\lambda_3(\bfalpha)$ to be the smallest
   positive root of the equation
   \beql{lambda3}
        |\det(J(\bfalpha))|-3n\cdot n!\cdot
        V(U+3Vx)^{n-1}\cdot x = 0.
   \eeql
    The following lemma implies the existence of $\lambda_{\jc}$:
    \bleml{lambdaintjc}
        If box $B$ contains a simple root $\bfalpha$ and
         $w_B< \lambda_3(\bfalpha)$ then  $\jc(B)$ succeeds.
    \eleml
    
    Thus we may choose $\lambda_{\jc}\as
    	\min_{\bfalpha\in \zero(2B_0)} \lambda_3(\bfalpha)$
    and set
	$$\ell_1\as \min\set{ \lambda_\jc, \lambda_\mk}$$

	\blemT{{{\bf Lemma A}}}{A}
		If $\#(B)>0$ and $w_B \le \ell_1$ then
		$\mk(B)$ and $\jc(B)$ holds.
	\elemT

	\bcor
		Each root in $B_0$ 
		will be output in a box of width $>3\ell_1/2$.
	\ecor

	Let $R_0\ib 2B_0$ be a region that excludes discs around roots:
		$$R_0\as 2B_0 \setminus
		\bigcup\nolimits_{\bfalpha\in \zero(2B_0)}
			\mathring{\Delta}(\bfalpha,\ell_1)$$
    where $\mathring{\Delta}$ is the interior of $\Delta$.
    Denote the zero set of $f_i$ as $\calS_i$ for $i=1\dd n$
	and define $d_0\as \inf_{p\in R_0} \max_{i=1}^n \sep(p,\calS_i)$.
    Since all the roots in $2B_0$ are removed from the set $R_0$,
    we can verify that $\max_{i=1}^n \sep(p,\calS_i)>0$ for all $p\in R_0$.
    Combining with the compactness of $R_0$, we obtain 
     $d_0>0$.
	Finally we set
		$$\lambda_\cz \as \frac{d_0}{2\sqrt{n}}.$$
    	
	\blemT{{{\bf Lemma B}}}{B}
		Suppose $\#(B)=0$.  If \\
        $\sep(\bfm_B,\zero(2B_0))\ge \ell_1$,
        \newOld{
        $\cz(B)$ succeeds when $w_B\le \lambda_\cz$. 
        If $\sep(\bfm_B,\zero(2B_0)) < \ell_1$,
        $\cz(B)$ succeeds when $w_B\le 
        \efrac{2}\min\{\lambda_\cz, \ell_1\}$. 
        }{if $w_B\le \lambda_\cz $ 
		then $\cz(B)$ holds.}
	\elemT

    The Lemma follows naturally from Lemma A and B:
	\blemT{{{\bf Lemma C}}}{C}
		Every box produced by the \miranda\ has width
		$\ge \efrac{4}\min\set{\lambda_\cz,\lambda_\jc,\lambda_\mk}$.
	\elemT

    \subsection{Sure Success for \protect\intcz$(B)$ and 
        \protect\intjc$(B)$}
We now consider the interval tests 
    $\intjc$ and $\intcz$ under the assumption that
    the underlying interval forms involved are Lipschitz.
    Let $\wh{L}$ be a global Lipschitz constant for 
    $\intbox f_i$ and $\intbox \frac{\partial f_i}{\partial x_j}$
    for all $i,j=1\dd n$ in $3B_0$.
    We will develop corresponding bounds $\lambda_{\intjc}$,
    $\lambda_{\intcz}$.
    Observe that if we replace the bounds 
    $\lambda_\mk$, $\lambda_\jc$, $\lambda_\cz$ in the abstract version
    by the bounds $\lambda_\intmk,\lambda_\intjc,\lambda_\intcz$,
    all the statements and proofs
    in the previous section remain valid. So
    in this section, we do not repeat the statements, 
    except to give the bounds $\lambda_\intjc$ and $\lambda_\intcz$.

    First look at the test $\intjc(B)$.
    With the same arguments as in abstract level, we obtain 
    	$$\lambda_{\intjc}\as
    		\min_{\bfalpha\in \zero(2B_0)} \lambda_4(\bfalpha)$$
   where $\lambda_4(\bfalpha)$ is the smallest positive root of the 
	   \beql{lambda4}
	        |\det(J(\bfalpha))|-3n\cdot n!\cdot
	        \wh{L}(U+3\wh{L}x)^{n-1}\cdot x = 0.
	   \eeql
    With $\lambda_\intjc$ and $\lambda_\intmk$, we get an 
    interval analogue of Lemma~A:
    \blemT{{\bf Lemma \intbox A}}{A'}
		If $\#(B)>0$ and $w_B \le \ell'_1$
        with $\ell'_1\as$ $\min\set{
        \lambda_\intjc, \lambda_\intmk}$, then
		$\intmk(B)$ and $\intjc(B)$ succeeds.
	\elemT
    Next look at the test $\intcz(B)$.
    Arguing as in the abstract level, we
    only consider the boxes in the region 
		$R_0'\as 2B_0 \setminus
		\bigcup\nolimits_{\bfalpha\in \zero(2B_0)}
			\mathring{\Delta}(\bfalpha,\ell'_1)$
        with $\ell'_1\as \min\set{\lambda_\intjc,\lambda_\intmk}$.
    Define $u\as \inf_{\bfp\in R'_0}\max_{i=1}^n\frac{|f_i(\bfp)|}{\wh{L}} $.
    It is easy to see that $\max_{i=1}^n\frac{|f_i(\bfp)|}{\wh{L}}>0$ 
    for any $\bfp\in R'_0$.
    Since the function $|f_i(\bfx)|$ is continuous and
    the set $R'_0$ is compact, we obtain that $u>0$.
    Setting $\lambda_\intcz\as \frac{u}{2}$, we have the following lemma:

    \blemT{{\bf Lemma~\intbox B}}{B'}
     \newOld{
		Suppose $\#(B)=0$.  
        Let \\
        $\ell'_1\as \min\set{\lambda_\intjc,\lambda_\intmk}$.
        If $\sep(\bfm_B,\zero(2B_0))\ge \ell'_1$,
        then $\intcz(B)$ succeeds when $w_B\le \lambda_\intcz$. 
        If $\sep(\bfm_B,\zero(2B_0)) < \ell'_1$,
        then $\intcz(B)$ succeeds when $w_B\le 
        \efrac{2}\min\{\lambda_\intcz, \ell'_1\}$. 
        }{Let $\sep(\bfm_B,\zero(2B_0))> \ell'_1$
        with $\ell'_1\as \min\set{\lambda_\intjc,\lambda_\intmk}$.
        If $\#(B)=0$ and
        $w_B\le \lambda_\intcz$, then  $\intcz(B)$ succeeds.}
    \elemT
    
    Combining Lemma~\intbox A and Lemma~\intbox B, we obtain:
    \blemT{{{\bf Lemma \intbox C}}}{C'}
	Every box produced by the \intbox\miranda\ has width
	$\ge \efrac{4}\min\set{\lambda_\intcz,\lambda_\intjc,\lambda_\intmk}$.
	\elemT    

\section{Effective \miranda}
	We now extend our results from $\intbox\miranda$
	to $\wtbox\miranda$ by introducing the effective tests
		$\wtmk(B)$, $\wtjc(B)$ and $\wtcz(B)$.
   \newOld{
    Recall that the difference between the effective version and
    the interval version is that the former uses dyadic numbers 
    instead of real numbers. In $\wtbox\miranda$, this difference 
    is reflected in 2 places: 
    (2) the preconditioning matrix $J\inv(\bfm)$ in 
    $\wtmk(B)$ is approximated by $\wt{J}\inv(\bfm)$ with dyadic entries 
    (1) each box form $\intbox h(B)$ (including $\boxm h(B)$) is 
    outwardly rounded to the effective form $\wtbox h(B)$
    whose endpoints are dyadic numbers.

    The main issue in $\wtbox\miranda$ is the accuracy of 
    $\wt{J}\inv(\bfm)$ and $\wtbox h(B)$.
    Based on the following requirements, we claim that similar lower
    bounds as in Section~\ref{complexity} for $\wtbox\miranda$ is 
    still achievable.
    We require
    {\small
        \beql{R1}\tag{R1}
        \|\wt{J}\inv(\bfm) - J\inv(\bfm) \|\le 
        \frac{1}{12n^2\cdot \|J(\bfm)\|}, 
        \eeql
        \beql{R2}\tag{R2}
        q(\intbox h(B),\wtbox h(B))\le \efrac{16}w_B.
        \eeql
    }
    }{Inside these tests are various box forms,
    say $\intbox h(B)$.
	Recall that they are actually mean value forms $\boxm h(B)$
    (we write ``$\intbox h(B)$'' for simplicity).
    We convert each $\intbox h(B)$ 
	to its effective version $\wtbox h(B)$, 
	which has dyadic endpoints and satisfies
	$\intbox h(B)\ib\wtbox h(B)$.  
    The main issue is the
	accuracy of the effective forms, which we express
	by upper bounds on the
	Hausdorff distance $q(\intbox h(B),\wtbox h(B))$.
    It is always bounded as a linear function of the width $w_B$,
	i.e., $q(\intbox h(B),\wtbox h(B))=O(w_B)$.
	However, we cannot stop here -- the implicit constant in
	the asymptotic notation must be made explicit
    for implementation purposes.
    }
    More precisely, for (R2) we require that 
    in $\wtcz(B)$,
		$q(\intbox f_i(B),$
			$\wtbox f_i(B)) \le$ $\efrac{16}w_B$
	for $i=1\dd n$;
    in $\wtmk(B)$, 
		$q(\intboxm g_i(2B^{\pm}_i),$
			$\wtboxm g_i(2B^{\pm}_i))$ $\le \efrac{16}w_B$
	for $i=1\dd n$;
	in $\wtjc(B)$, 
		$q(\intbox J_{ij}(3B),$
			$\wtbox J_{ij}(3B))\le$ $\efrac{16}\cdot 3w_B$
	for each entry $\wtbox J_{ij}(3B)$ of $\wtbox J(3B)$. 
    
    \newOld{It is easy to verify that,
    for $\wtjc$ and $\wtcz$, we can get similar bound as 
    $\lambda_{\intjc}$ and $\lambda_{\intcz}$ by simply replacing 
    $\wh{L}$ with $\wh{L}+\efrac{8}$.
    For $\wtmk$, similar bound as $\lambda_{\intmk}$ can also be
    obtained by modified the definition of $\lambda_1(\bfalpha)$.
    }
    {We get effective versions of all our lemmas and theorems,
	with modified constants such as
    $\lambda_{\intjc}$ and $\lambda_{\intcz}$.}
	\ifNeedSpace{
	It is straightforward to confirm that $\lambda_{\mk}$ still holds
	for $\wtmk$.
	As for $\wtjc$ and $\wtcz$,
	it suffices to increase the Lipschitz constant
	by a factor $\efrac{8}$, and the proofs for
	$\lambda_{\intjc}$ and $\lambda_{\intcz}$ will go through.
	}
	\ignore{
	    TO DO: how to treat
	    the inverse Jacobian $J\inv$  in MK test
	}%

\section{Conclusion}
	We have provided the first effective subdivision
	algorithm \miranda\ for isolating simple real roots of a system
	of equations $\bff=\0$, provided $\bff$ and its derivatives have
	interval forms.  Our result are novel for its completeness
	(previous algorithms need $\eps$-termination and has no isolation
	guarantees), its generality (going beyond the polynomial case),
	and its complexity analysis (going beyond
	termination proofs).  We also contributed to
	the theory of subdivision algorithms by formalizing a 3-level
	description to provide a pathway from abstract algorithms
	to effective ones.  Given that many existing numerical algorithms
	still lack effective versions, this is a promising line of work.
	In the future, we plan to implement and develop our
	algorithm into a practical tool.  

\section*{Acknowledgments}
The first author wishes to thank her supervisor Professor Dongming Wang
for his support and encouragement. 
The second author is thankful for hospitality 
and support during a sabbatical leave:
Professor Jing Yang
(Guangxi University of Nationalities, Nanning),
Drs.~Xiaoshan Gao and Jinsan Cheng
(Chinese Academy of Sciences, Beijing) and
Professors Dongming Wang and Chenqi Mou
(Beihang University).
The following support is acknowledged:
NSF Grants \#CCF-1423228 and \#CCF-1564132;
a Chinese Academy of Science President's International Fellowship
Initiative (2018), and
Beihang International Visiting Professor Program No. Z2018060.

\bibliographystyle{plain}

\newpage
\section{Appendix: All Proofs}

    \bgenDIY{{\sc Lemma 4.1 (MVT).}}
    {Given two points $x,y\in \RR^n$, we have:\\
	(a)
	\beql{MVT1}
        J(\bfx) = J(\bfy) \pm \K{[\bfx,\bfy]}\|\bfx-\bfy\|
	\eeql
	(b)
	\beql{MVT2}
        \bff(\bfx) - \bff(\bfy) = (J(\bfy) \pm \K{[\bfx,\bfy]}
			\|\bfx-\bfy\|)\cdot(\bfx-\bfy)
	\eeql
    }
	\bpf
	(a) 
	We apply the Mean Value Theorem to each entry
	$J_{ij}= \frac{\partial f_i}{\partial x_j}$:
	\begin{equation*}
	\begin{aligned}
	J_{ij}(\bfx) &= J_{ij}(\bfy)+
	           \nabla J_{ij}(\wt{\bfy})\cdot(\bfx - \bfy)
	           \quad\quad \text{with}\ \wt{\bfy}\in [\bfx, \bfy]\\
	             &= J_{ij}(\bfy) \pm K([\bfx,\bfy])_{ij}\|\bfx-\bfy\|
	\end{aligned}
	\end{equation*}
	(b) We apply the Mean Value Theorem twice.  The first application
	gives:
	\begin{equation*}
	\begin{aligned}
		f_i(\bfx)-f_i(\bfy)
	    	&= \nabla f_i(\wt{\bfy})\cdot(\bfx - \bfy) \\
        &= (J_{i1}(\wt{\bfy})\dd J_{in}(\wt{\bfy}))\cdot(\bfx-\bfy)
	\end{aligned}
	\end{equation*}
	where $\wt{\bfy}\in [\bfx,\bfy]$ and
	$J_{ij}\as \frac{\partial f_i}{\partial x_j}$.
	Applying the Mean Value Theorem again to each $J_{ij}(\wt{\bfy})$:
	\begin{equation*}
	\begin{aligned}
	J_{ij}(\wt{\bfy})
	    	&= J_{ij}(\bfy)+
	            \nabla J_{ij}(\wh{\bfy})\cdot(\bfy - \wt{\bfy})
	            \quad\quad \text{with}\ \wh{\bfy}\in [\bfy, \wt{\bfy}]\\
		&= J_{ij}(\bfy) \pm K([\bfx,\bfy])_{ij}\|\bfx-\bfy\|\\
	\end{aligned}
	\end{equation*}
	Hence
	\begin{equation*}
	\begin{aligned}
	f_i(\bfx)-f_i(\bfy)
        = (& J_{i1}(\bfy)
	    \pm K([\bfx,\bfy])_{i1}\|\bfx-\bfy\|\dd \\
           & J_{in}(\bfy) \pm
	    K([\bfx,\bfy])_{in}\|\bfx-\bfy\|)\cdot (\bfx-\bfy)
	\end{aligned}
	\end{equation*}
	for $i=1\dd n$.
	This proves \refeq{MVT2}.
	\epf
    
    \bgenDIY{{\sc Lemma 4.2.}}
    {For any simple root $\bfalpha$ of $\bff$, $\lambda_1(\bfalpha)$
	is well-defined.
    }
	\bpf
	Note that $s(0)$ is well-defined since $\bfalpha$ is a simple root.
	We also deduce that $s(0)<0$ and that
	$s(r)=s(0)$ for all $r<0$.  Thus $\lambda_1(\bfalpha)>0$
	if it is well-defined.
    Let $r\aast$ be the smallest radius such that $\Delta(\bfalpha,r\aast)$
    contains a critical point; if $\bff$ has no critical point, then
    $r\aast$ is defined to be $\infty$.
    It follows that $s(r\aast)=r\aast - \frac{1}{\infty} = r\aast$.
    Thus $s(0)<0<s(r\aast)$.
    From the fact that
    $\jkd{\bfalpha,2\sqrt{n}r}{\bfalpha,2\sqrt{n}r}$
    is a continuous non-decreasing
	function of $r$ in the range $[0,r\aast)$, we conclude
	that there exists some $r\in (0,r\aast)$ such that $s(r)=0$.
    	\epf

    \bgenDIY{{\sc Lemma 4.3.}}
        {Let $B$ be a box containing a simple root $\bfalpha$ of $\bff$
        and $\bfm\in B$ with $J\inv(\bfm)$ well-defined.
	    If $w_B\le \lambda_1(\bfalpha)$, the preconditioned system 
        $\bfg_B \as J\inv(\bfm)\bff = (g_1\dd g_n)$ 
	    satisfies that for all $i=1\dd n$,
	    	$$g_i(2B^+_i) \ge \frac{w_B}{4},\qquad
	    	g_i(2B_i^-) \le -\frac{w_B}{4}.$$
        }
	\bpf
	    Let $\bfx$ be a point on the boundary of the box $2B$.
	    Then
	    {\scriptsize\begin{align*}
        & \bfg_B(\bfx)\\
		&=J\inv(\bfm)\bff(\bfx)
            & \hspace{-80pt}\textrm{(by definition of $\bfg_B$)}\\
		&=J\inv(\bfm)(\bff(\bfalpha) + (J(\bfalpha) \pm
            \K{[\bfx,\bfalpha]}\|\bfx - \bfalpha\|) \cdot (\bfx - \bfalpha))
			& \hspace{-80pt}\textrm{(by MVT \refeq{MVT2})}\\
		&=J\inv(\bfm)(J(\bfalpha) \pm
		\K{[\bfx,\bfalpha]}\|\bfx - \bfalpha\|) \cdot (\bfx - \bfalpha)
			& \hspace{-80pt}\textrm{(since $\bfalpha$ is a root)}\\
        &=J\inv(\bfm)(J(\bfm)\pm \K{[\bfalpha,\bfm]}\|\bfalpha-\bfm\| \pm
		\K{[\bfx,\bfalpha]}\|\bfx - \bfalpha\|) \cdot (\bfx - \bfalpha)\\
        &	& \hspace{-80pt}\textrm{(by MVT \refeq{MVT1})}\\
		&= J\inv(\bfm)(J(\bfm) \pm 3K(2B)\|\bfx - \bfalpha\|)
                    \cdot (\bfx - \bfalpha)
			& \hspace{-80pt}\textrm{(since $\|\bfm-\bfalpha\|
			\le2\|\bfalpha-\bfx\|$)}\\
            &= (\1 \pm 3J\inv(\bfm)K(2B)\|\bfx - \bfalpha\|)
			\cdot (\bfx - \bfalpha)
			& \hspace{-80pt}\textrm{($\1$ is the identity matrix)}.
	    \end{align*}
	    }
	
	   The $i$-th component in $\bfg_B(\bfx)$ is the $g_i$; thus
	    \begin{align*}
	  g_i(\bfx) &=  (x_i-\alpha_i)
			\pm 3(J\inv(\bfm)K(2B)\|\bfx - \bfalpha\|)
			\cdot (\bfx - \bfalpha).
	    \end{align*}
	In the following, we write $\lambda_1$ for $\lambda_1(\bfalpha)$
	and note that
	$\bfalpha\in B$ and $w_B\le \lambda_1$ implies
	\beql{m2b}
		\grouping{
			m&\in& \discalpha\\
			2B&\ib& \discalpha
			}.
	\eeql
	    Thus:
	    {\footnotesize
	    \begin{flalign*}
        &\Big|g_i(\bfx) - (x_i - \alpha_i)\Big| \\
		&\le 3\|J\inv(\bfm)K(2B)\|\cdot
                    \|\bfx - \bfalpha\|\sum^n_{j=1}|x_j - \bfalpha_j| \\
        &\le \frac{9}{2}nw_B\|J\inv(\bfm)K(2B)\|\cdot \|\bfx - \bfalpha\|
            & \hspace{-16pt}\textrm{(as $\sum^n_{j=1}|x_j - \bfalpha_j|\le
		    \frac{3}{2}nw_B$)}\\
		& \le \frac{27}{4}nw_B^2  \|J\inv(\bfm)K(2B)\|
		    & \textrm{(as $\|\bfx - \bfalpha\| \le \frac{3}{2}w_B$)}\\
		& \le \frac{w_B^2}{4}\Big(27n
			\|J\inv(\discalpha)
			 K(\discalpha)\|\Big)
		    & \textrm{from \refeq{m2b}}\\
		& = \frac{w_B^2}{4} \cdot \efrac{\lambda_1}
		    & \textrm{(definition of $\lambda_1$)}\\
		& \le \frac{w_B}{4}
		    & \textrm{(since $w_B\le \lambda_1$)}.
		    \label{eq:wb4}
            \end{flalign*}
	    }
	This last inequality gives
		\beql{wb4}
		\Big|g_i(\bfx)- (x_i-\alpha_i)\Big|\le w_B/4.
		\eeql
	It remains to show that $g_i(2B^+_i)\ge \frac{w_B}{4}$
	(the proof that $g_i(2B^-_i)\le -\frac{w_B}{4}$ is similar).
	This amounts to proving $g_i(\bfx)\ge \frac{w_B}{4}$
	holds for all $\bfx\in 2B^+_i$.
	First we note that
		\beql{xi} x_i-\alpha_i \ge w_B/2	\eeql
	since
	$\bfx\in 2B^+_i$ and $\bfalpha\in B$.
	The inequalities \refeq{wb4} and \refeq{xi}
	together implies $g_i(\bfx)$ and $x_i-\alpha_i$
	must have the same sign.  
	Since $x_i-\alpha_i$ is positive, we conclude
	that $g_i(\bfx)$ must be positive.
	Combined with \refeq{wb4} and \refeq{xi}, we
	conclude that $g_i(\bfx)\ge w_B/4$, as claimed.
	\epf

    \bgenDIY{{\sc Lemma 4.4.}}
    {Let $f$ be a continuously differentiable function defined on a 
    convex region 
    $S\ib \RR^n$.  
        Then 
    $\sum_{k=1}^n \Big|\intbox \frac{\partial f}{\partial x_j}(S)\Big|$
    is a Lipschitz constant for $\boxm f$ on $S$.
    }
    \bpf
    Recall that
    $\boxm f(B) = f(\bfm(B))+\intbox\nabla f(B)^T\cdot (B-\bfm(B))
    = f(\bfm(B))+\efrac{2}w_B\cdot \sum_{k=1}^n\intbox 
    \frac{\partial f}{\partial x_j}(B)$
    for any $B\ib S$.
    Thus $w(\boxm f(B)) = \efrac{2}w_B\cdot w(\sum_{k=1}^n
    \intbox \frac{\partial f}{\partial x_j}(B))
     =\efrac{2}w_B\cdot\sum_{k=1}^n
    w(\intbox \frac{\partial f}{\partial x_j}(B))
    \le w_B\cdot \sum_{k=1}^n\Big|\intbox\frac{\partial f}{\partial x_j}(B)\Big|
    \le w_B\cdot 
    \\\sum_{k=1}^n\Big|\intbox\frac{\partial f}{\partial x_j}(S)\Big|$. The lemma follows.
    \epf

    \bgenDIY{{\sc Theorem 4.6.}}
        {Let $B$ be a box containing a simple root $\bfalpha$
        of width $w_B\le \lambda_{1}(\bfalpha)$ and 
        $\bfm\in B$ with $J\inv(\bfm)$ well-defined.
        \benum[(a)]
        \item 
	    If $w(\intbox\frac{\partial g_i(2B_i^+)}{\partial x_j})
			\le \frac{1}{32n}$
	    for each $j=1\dd n$, then $\intmk(B)$ succeeds
	    with $\bfg_B := J\inv (\bfm)\bff$.
        \item
        	If $w_B\le \lambda_2(\bfalpha)$ with
        	$\lambda_2(\bfalpha)\as\min\set{\lambda_1(\bfalpha),
        		\wh{\lambda}_1(\bfalpha)}$,
        	then $\intmk(B)$ succeeds.
        \eenum
        }
	\bpf
        (a) We show the first part of the theorem.
        In \refLem{theoremMK}, it is proven that when
        $w_B\le \lambda_1(\bfalpha)$,
	    it holds 
	        $g_i(B_i^+)\ge\frac{w_B}{4}$ and $g_i(B_i^-)\le-\frac{w_B}{4}$.
	    From \refPro{quadConverg}, we have
	    \begin{align*}
	        q(\boxm {g_i}(2B_i^+), g_i(2B_i^+))
		&\le 2 w(2B)\sum_{j=1,j\neq i}^n w(\intbox\frac
                    {\partial g_i(2B_i^+)}{\partial x_j}) \\
	        &\le 4nw_B\cdot\max_{1\le j\le n,j\neq i}
                w(\intbox\frac{\partial g_i(2B_i^+)}{\partial x_j}).
	    \end{align*}
	    By the convergence property of box functions,
        $w(\intbox\frac{\partial g_i(2B_i^+)}{\partial x_j})$
        approaches $0$ when $w_B$ approaches $0$
	    for $j=1\dd n$.
	    Thus when $w_B$ is small enough, we have
	      $w(\intbox\frac{\partial g_i(2B_i^+)}{\partial x_j})
	      		\le \frac{1}{32n},
	    \forall j=1\dd n$. 
        Then
	     \begin{align*}
	        \boxm {g_i}(2B_i^+)
	         &\ge g_i(2B_i^+)-q(\boxm {g_i}(2B_i^+), g_i(2B_i^+)) \\
	         &\ge\frac{w_B}{4}-4nw_B\cdot\frac{1}{32n} = \frac{w_B}{8}>0.
	     \end{align*}
	   Similar argument applies to $\boxm {g_i}(2B_i^-)$.
       This gives the first part of the theorem.

       (b) Now we prove the second part of the theorem.
       From the proof of the first part, it suffices to prove that
       when $w_B\le \lambda_2(\bfalpha)$,
       the inequality
	   $w(\intbox\frac{\partial g_i(2B_i^+)}{\partial x_j})
            \le \frac{1}{32n}$
            holds for all $i,j=1\dd n$.
       To show this, we observe that
       {\footnotesize
       \begin{align*}
       & w(\intbox\frac{\partial g_i(2B_i^+)}{\partial x_j}) \\
       & = \sum_k [J\inv(\bfm)]_{ik}\cdot \intbox \frac{\partial f_j}{
           \partial x_k}(2B_i^+) \\
       & &\hspace{-50pt}\textrm{($[J\inv(\bfm)]_{ik}$ are the entries of $J\inv(\bfm)$)}\\
       & < \sum_k \|J\inv(\discalpha)\|\cdot \intbox \frac{\partial f_j}{
           \partial x_k}(2B_i^+)
       & \textrm{($2B\subset\discalpha$)} \\
       & <\|J\inv(\discalpha)\|\cdot 2nL w_B
       & \textrm{($\intbox \frac{\partial f_j}{\partial x_k}$ are Lipschitz on
        $\discalpha$)}\\
       & \le \frac{1}{32n}
       &\textrm{($w_B\le \efrac{64n^2L\cdot
       \|J\inv(\discalpha)\|}$)}.
       \end{align*}
        }
	   \epf

    \bgenDIY{{\sc Theorem 5.2.}}
        { If both $\jcs(B)$ and
		$\mk(\frac{3}{2}B)$ succeed then $\#_\bff(3B)=1$.
        }
    \bpf
        From \cite{franek-ratschan:topotest:12}, the success of
        $\mk(\frac{3}{2}B)$ implies 
            $$\sum_{\bfy\in\zero(3B)} {\rm sign}( \det J_\bff(\bfy))
                 = \pm1$$
        where ${\rm sign}( \det J_\bff(\bfy)))$ is the sign of 
        $\det J_\bff(\bfy)$.
        By the success of $\jcs(B)$, we further know that 
        ${\rm sign}(\det J_\bff(\bfy))$ is the same for all $\bfy \in 3B$.
        Thus there is only one root in $3B$.
    \epf

    \ignore{
    \bgenDIY{{\sc Lemma XXX.}} 
        {If $\jc(B)$ succeeds, $3B$ contains at most one root.}
    \bpf
        Prove by contradiction. Suppose that there are $2$ roots 
        $\bfx, \bfy$ in $3B$.
        Then for each $i\in [1,n]$, there exist a point 
        $\bfz_i \in [\bfx, \bfy]$ such that 
        $\nabla f_i(\bfz_i) \cdot (\bfx - \bfy) = 0$.
        Thus the $n$ rows of the matrix 
        $[\nabla f_1\dd \nabla f_n]^T$ cannot be independent.
        It follows $\det [\nabla f_1\dd \nabla f_n]^T = 0$.
        And $[\nabla f_1\dd \nabla f_n]^T \in J_\bff(3B)$ since 
        $\bfz_i\in 3B$ for all $i\in [1,n]$.
        Thus  $\jc(B)$ fails.
    \epf
	}%

    \bgenDIY{{\sc Theorem 6.1 (Partial Correctness).}}
    {If \miranda\ halts, the output queue $P$ is correct.
	   }
	\bpf
	\ignore{
    	    First we prove the termination of the algorithm. We prove this
	    by contradiction.
	    Suppose the algorithm does not terminate, then we can find an
	    infinite series of boxes $B_1,B_2,\cdots$ such that
	    $B_1\supset B_2\supset \cdots$.
	    Since width of the boxes is halved each time,
	    this infinite series must converge to a point, denoted as
	    $p$.
	    We divide the proof into 2 cases depending on whether
	    or not $p$ is a root of $\bff$.
	
	    If $p$ is a root of $\bff$, then from \refThm{MK}, the test
	    \intmk\ will succeed on any sufficient box containing $p$.
	    That is, there exists an integer $A_1$ such that
            \intmk$(B_i)$ succeeds for $\forall i\ge A_1$.
	    And by assumption,  $\det(J(p))\neq 0$.
	    Hence the convergence property of box functions implies
	    that there exists an integer $A_2$ such that $\jc(B_i)$ succeeds
	    for $\forall i\ge A_2$.
            Therefore, if $i=\max\{A_1,A_2\}$, $B_i$ would satisfy
            the $\mk$ and $\jc$ tests and be output as an isolating box.
            This is a contradiction.
	
	    If $p$ is not a root of $\bff$, then by the convergence property of
            box functions,
	    there exists an integer $A_3$ such that $\cz(B_i)$ succeeds
	    for $\forall i\ge A_3$.
	    That is, the box $B_i$ will be discarded for $\forall i\ge A_3$.
	    This is also a contradiction.
	}%
    \ignore{ 
	    We first argue that the partial correctness of
	    \intbox\miranda\ and $\wtbox$\miranda\ follows 
	    from the partial correctness of \miranda\ by
	    the general observation\footnote{
		If we were proving termination,
		the reverse implication hold: if $\wtbox$\miranda\
		terminates than \miranda\ terminates.
	    }
	    that the predicates in \miranda\ are one-sided,
	    and (as can be verified below) none of our arguments
	    are predicated upon the {\em failure} of the tests.
	    We need to further note that for the effective version,
	    we must assume that the ROI $B_0$ is a dyadic box,
	    so that all subdivisions are done without approximation.
    }%
        Firstly, we note that each output box in $P$ is isolating.
	    A box $2B$ is output in line 11 upon
	    passing $\mk(B)$. This is inside the inner while loop for subboxes
	    of some $B'$ which passes $\jc(B')$.  But $\mk(B)$ implies
	    $\#(2B)\ge 1$ and $\jc(B')$ implies $\#(3B)\le 1$.
	    Thus $\#(2B)=1$.
	
	    Next we claim no root is output twice in $P$.  
	    {\em This follows by showing that if $2B$ and $2B'$ are
	    output, then their interiors are disjoint.}
	    It does not matter
	    if the boundaries of $2B$ and $2B'$ intersect because
	    there are no roots on their boundary -- this is
	    ensured by the success of the Simple Miranda test on
	    these output boxes.
	    The reason for our concern comes from the fact that,
	    although the boxes in $Q$ have pairwise disjoint interiors,
	    each $B$ in $Q$ can cause a larger box $2B$ to be output.
	
	    {\em CLAIM: Suppose $2B$ is output in line 11.
	    Then immediately after line 12, every box
	    $B'$ in $Q$, the interior of $2B'$ is disjoint from $2B$.
	    }
	    Pf: Suppose the interior of $2B'$ intersects $2B$.
	    By the priority queue property, we have
	    $w(B')\le w(B)$.  It follows that $B'$ actually
	    is contained in the annulus $3B\setminus B$.
	    This follows from two facts\footnote{
		Here \dt{aligned boxes} means those that can
	    	arise by repeated subdivision of $B_0$.
		Clearly $B$ and $B'$ are aligned, but $kB$ and $kB'$ are
		not aligned for any $k>1$.
	    } about aligned boxes: (a) any two aligned boxes
	    have disjoint interiors or have a containment relationship,
	    and (b) $3B\setminus B$  is a union of $8$ aligned boxes.
	    If $B'$ is contained in this annulus, then
	    line 12 would have removed it.
	    This proves our claim.
	
	    Finally we must show that
		 $$\calZ_\bff(B_0)
		 	\stackrel{(*)}{\subseteq} \bigcup_{B\in P}\calZ_\bff(B)
		 	\stackrel{(**)}{\subseteq} \calZ_\bff(2B_0).$$
	    The second containment (**) is immediate because all
	    our output boxes have the form $2B$ where $B$ is an aligned
	    box.  Such boxes are contained in $2B_0$.
	    To show (*),
	    it suffices to prove that if $B'$ is a discarded box,
	    then either $B'$ has no roots, or any
	    root in $B'$ is already output.
	    From the algorithm, a box $B'$ is discarded in two lines:
	    The first is Line 4, when $\cz(B')$ succeeds.
	    But this implies $B$ has no roots.  
	    The second is Line 12 of the algorithm.
	    Since $B'$ in contained in $3B$ (where $2B$ is the output).
	    We know that $\jc(B)$ holds, and
	    thus there is at most one root in $3B$.
	    So if $B'$ contains any root, it must be the root already
	    identified by $2B$.  Thus, all discarded boxes are justified.
		\ignore{
	    It remains to prove that any two distinct boxes in $P$
	    do not contain a same root of $\bff$.
	    To this end, we will show that the interior of any two boxes in $P$
	    do not overlap.
	    Take $B_1',B_2'\in P$ and without loss of generality,
	    suppose that $B_1'$ is output earlier than $B_2'$.
	    By the algorithm there exists $B_1$ such that
	    $B_1'\subset 2B_1$ and \intjc$(B_1)$ succeeds.
	    Since all the boxes in $3B_1\setminus B_1'$ are discarded,
	    the distance from $B_1'$ to any remained box
	    is at least $\frac{1}{2}w(B_1)$.
	    Hence the distance from $B_1'$ to $\frac{1}{2}B_2'$ is at least
	    $\frac{1}{2}w(B_1)$.
	    Now it suffices to show that
	      $\frac{1}{4}w(B_2')\le\frac{1}{2}w(B_1)$.
	    This is true because $B_1$ is a box of biggest size in $Q$,
	    thus $w(B_1)\ge w(\frac{1}{2}B_2')$.
	    }%
	\epf

    \bgenDIY{{\sc Lemma 7.1.}}
        { If box $B$ contains a simple root $\bfalpha$ and
         $w_B< \lambda_3(\bfalpha)$ then  $\jc(B)$ succeeds.
        }
    \bpf
        The fact $\bfalpha\in B$ implies
            $J(\bfalpha)\in J(3B)$.
        Since $\bfalpha$ is a simple root, we have
        $\det(J(\bfalpha))\neq 0$, and thus
        $\lambda_3(\bfalpha)\neq 0$.
        From the definition of $\lambda_3(\bfalpha)$,
        we know that if $w_B<\lambda_3(\bfalpha)$,
        then $3n\cdot n!\cdot V(U+3Vw_B)^{n-1}\cdot w_B
        <|\det(J(\bfalpha))|$, and thus
        $w(\det(J(3B)))<|\det(J(\bfalpha))| $.
        It follows $0\notin \det(J(3B))$.
        The test $\jc(B)$ succeeds.
    \epf

    \bgenDIY{{\sc Lemma 7.4} ({\bf Lemma B}).}
        {
		Suppose $\#(B)=0$.  If 
        $\sep(\bfm_B,\zero(2B_0))$ $\ge \ell_1$,
        $\cz(B)$ succeeds when $w_B\le \lambda_\cz$. 
        If $\sep(\bfm_B,\zero(2B_0)) < \ell_1$,
        $\cz(B)$ succeeds when $w_B\le 
        \efrac{2}\min\{\lambda_\cz, \ell_1\}$. 
        }
	\bpf
	    First consider the case $\sep(\bfm_B,\zero(2B_0))\ge\ell_1$.
        In this case, we have
        $\bfm_B\in R_0$. Thus $\max_{i=1}^n \sep(\bfm_B,\calS_i)\ge d_0$.
        Combining $w_B\le \lambda_\cz$, it follows
        that $\max_{i=1}^n \sep(B,\calS_i)\ge d_0- \sqrt{n}w_B
        \ge d_0/2$. Hence there exists $i\in [1,n]$ such that
        $B\cap \calS_i = \emptyset$.
        Thus $\cz(B)$ holds.

	    Then consider the case $\sep(\bfm_B,\zero(2B_0)) < \ell_1$.
        Without loss of generality, assume that 
         $\sep(\bfm_B,\bfalpha) < \ell_1$.
        Since $w_B\le \half\min\set{\lambda_\cz,\ell_1}$,
        by the Corollary to \refLem{A},
        we know that some box $B'$ containing $\bfalpha$
        of width $>\ell_1/2$ has been output (this uses
        the fact that we process the boxes in a
        breadth-first-manner).
        This output also removes all the boxes in the 
        process queue   
        that intersect the interior of $3B'$. 
        We can see that $B$ intersects $3B'$ because
        $w_{B'}>\ell_1/2$, $\bfalpha\in B'$ and 
         $\sep(\bfm_B,\bfalpha) < \ell_1$.
        Thus $B$ should have been removed. 
        Contradiction.
	\epf

    \ignore{
    \bgenDIY{{\sc Lemma 15} ({\bf Lemma C}).}
        {Every box produced by the \miranda\ has width
		$\ge \efrac{4}\min\set{\lambda_\cz,\lambda_\jc,\lambda_\mk}$.
        }
	\bpf
        We first give an equivalent statement of the lemma:
        for any box $B$ produced in the algorithm,
        if $w_B\le \half\min\set{\lambda_\cz,\ell_1}$,
        then $B$ has been output or rejected.
        In what follows, we will prove this equivalent statement.


		Case 1, $\#(B)>0$: Thus some ancestor of $B$
		should have output according to \refLem{A}, a contradiction.

		Case 2, $\#(B)=0$, $\sep(\bfm_B,\zero(2B_0))\ge \ell_1$.
			We contradict \refLem{B} directly.

		Case 3, $\#(B)=0$, $\sep(\bfm_B,\zero(2B_0)) < \ell_1$.
            Without loss of generality, assume that 
             $\sep(\bfm_B,\bfalpha) < \ell_1$.
			Since $w_B\le \half\min\set{\lambda_\cz,\ell_1}$,
            by the Corollary to \refLem{A},
			we know that some box $B'$ containing $\bfalpha$
			of width $>\ell_1/2$ has been output (this uses
			the fact that we process the boxes in a
			breadth-first-manner).
			This output also removes all the boxes in the 
            process queue   
			that intersect the interior of $3B'$. 
			We can see that $B$ intersects $3B'$ because
            $w_{B'}>\ell_1/2$, $\bfalpha\in B'$ and 
             $\sep(\bfm_B,\bfalpha) < \ell_1$.
            Thus $B$ should have been removed. 
			Contradiction.
	\epf
    }

    \bgenDIY{{\sc Lemma 7.7.} ({\bf Lemma \intbox B}).}
        {Let $\sep(\bfm_B,\zero(2B_0))> \ell'_1$
        with $\ell'_1\as$ $\min\set{\lambda_\intjc,\lambda_\intmk}$.
        If $\#(B)=0$ and
	$w_B\le \lambda_\intcz$, then  $\intcz(B)$ succeeds.
        }
    \bpf
    The Proof is similar to that of \refLem{B}.
        Since $\sep(\bfm_B,\zero(2B_0))>\ell'_1$, we have
        $\bfm_B\in R'_0$. 
        
        By the definition of $u$,
        we see $\max_{i=1}^n \frac{|f_i(\bfm_B)|}{\wh{L}}\ge u$,
        which means that $\exists j\in [1,n]$ such that 
        $\frac{|f_j(\bfm_B)|}{\wh{L}}\ge u$.
        By the inclusion property of box functions,
        $f_j(\bfm_B)\in \intbox f_j(B)$.
        Since $w_B\le \lambda_\intcz =\frac{u}{2}$, 
        we have $w(\intbox f_j(B))\le \wh{L}\cdot w_B\le \frac{u}{2}$.
        It follows
        that $\intbox f_j(B)\ge f_j(\bfm_B)-w(\intbox f_j(B))\ge
            u - \frac{u}{2}> 0$.
        Thus $0\notin \intbox f_j(B)$ and $\intcz$ holds.
    \epf
\end{document}

	\cite{ratschek-rokne:range:bk}
	\cite{neumaier:equations:bk-90}
	\cite{kearfott:bk}
	\cite{hentenryck+2:solving:97}
	\cite{hansen:interval:06}
	\cite{aberth:computable-analysis:bk-80}
	\cite{aberth:precise:bk} 

	\cite{mourrain-pavone:subdiv-polsys:09} 
	\cite{elber-kim:solver:01} 
	\cite{sherbrooke-patrikalakis:93} 
	\cite{garloff-smith:subdiv-algo:01}
	\cite{garloff-smith:bernstein-systems:01}
	\cite{moore-kioustelidis:test:80}
	\cite{dickenstein-emiris:poly-sys:05}
	
	REAL:		\cite{sagraloff-mehlhorn:real-roots:16}
			\cite{kobel-rouillier-sagraloff:for-real:16}
	COMPLEX:	\cite{becker+3:cisolate:18}
			\cite{becker+4:cluster:16}
			\cite{imbach-pan-yap:ccluster:18}
	\cite{lien-sharma-vegter-yap:arrange:14z} 
	\cite{bennett-papadopoulou-yap:minimization:16}
	\cite{yap-sagraloff-sharma:cluster:13}

	\cite{frommer2004comparison}
	\cite{frommer-lang:miranda-type-tests:05}
	\cite{goldsztejn2007comparison}